\documentclass[onefignum,onetabnum]{siamart171218}

\usepackage{amsmath}
\usepackage{multirow}
\usepackage{booktabs} 
\usepackage{bm}
\usepackage{lipsum}
\usepackage{amsfonts}
\usepackage{graphicx}
\usepackage{epstopdf}
\usepackage{float}
\usepackage{algorithmic}
\usepackage{subcaption}
\usepackage{geometry}
\allowdisplaybreaks
\usepackage{adjustbox}
\usepackage{hyperref}
\usepackage{cleveref}
\ifpdf
  \DeclareGraphicsExtensions{.eps,.pdf,.png,.jpg}
\else
  \DeclareGraphicsExtensions{.eps}
\fi

\newsiamremark{remark}{Remark}
\newsiamremark{hypothesis}{Hypothesis}
\crefname{hypothesis}{Hypothesis}{Hypotheses}
\newsiamthm{claim}{Claim}

\headers{IDR($s$) Method for Matrix Function Bilinear Forms}{Q. Xue, X. Yue, X.-M. Gu}

\title{Error Analysis of Krylov Subspace approximation Based on IDR($s$) Method for Matrix Function Bilinear Forms}

\author{
Qianqian Xue\thanks{School of Mathematics and Computational Science, Xiangtan University, Xiangtan, Hunan 411105, P.R. China and School of Mathematics, Southwestern University of Finance and Economics, Chengdu, Sichuan 611130, P.R. China(\email{222070100008@smail.swufe.edu.cn}).}
\and 
Xiaoqiang Yue\thanks{National Center for Applied Mathematics in Hunan, Key Laboratory of Intelligent Computing \& Information Processing of Ministry of Education, Hunan Key Laboratory for Computation and Simulation in Science and Engineering, Xiangtan University, Xiangtan, Hunan 411105, P.R. China
(\email{yuexq@xtu.edu.cn}).}
\and Xian-Ming Gu\footnotemark[4]\thanks{Corresponding author. School of Mathematics, Southwestern University of Finance and Economics, Chengdu, Sichuan 611130, P.R. China
(\email{guxianming@live.cn}, \email{guxm@swufe.edu.cn}).}}

\usepackage{amsopn}

\makeatletter
\newcommand*{\addFileDependency}[1]{
  \typeout{(#1)}
  \@addtofilelist{#1}
  \IfFileExists{#1}{}{\typeout{No file #1.}}
}
\makeatother

\ifpdf
\hypersetup{
  pdftitle={Error Analysis of Krylov Subspace approximation Based on IDR($s$) Method for Matrix Function Bilinear Forms},
  pdfauthor={Q. Xue}
  }
\fi

\begin{document}

\maketitle
%
\begin{abstract}
The matrix function bilinear form, namely $\mathbf{u}^\top f(A)\mathbf{v}$, appears in many scientific computing problems, where $\mathbf{u}, \mathbf{v} \in \mathbb{R}^n$, $A \in \mathbb{R}^{n\times n}$, and $f(z)$ is a given analytic function. The Induced Dimension Reduction (IDR($s$)) method was originally proposed to solve a large system of linear equations, and effectively reduces the complexity and storage requirement by dimensionality reduction while maintaining the numerical stability of the algorithm. In fact, the IDR($s$) method can generate an interesting Hessenberg decomposition. We make use of this decomposition to establish the numerical algorithm and a practical error estimation framework for the matrix function bilinear form. Based on an error analysis of the IDR($s$) approximation, the corresponding error expansion is derived. 
Crucially, we prove that, under mild conditions, the remainder term in this expansion converges to zero as the truncation order increases, thereby validating the error expansion. The leading computable term is then used as a practical a posteriori error indicator for general analytic functions. We present numerical experiments to support our theoretical findings and illustrate the efficacy of our proposed method, along with its stopping criterion, over traditional Bi-Lanczos and Arnoldi-based algorithms.
\end{abstract}



\section{Introduction}
This paper presents a reliable numerical method for computing the matrix function bilinear form, i.e.,
\begin{equation}
\mathbf{u}^{\top} f(A)\mathbf{v} ,
\label{1.1} 
\end{equation}
with the matrix \(A\in \mathbb{R}^{n\times n}\) being a large (and often sparse) matrix, and $\mathbf{u}$, $\mathbf{v}$ satisfying \(\|\mathbf{u}\|=1\), \(\|\mathbf{v}\|=1\), where \(\|\cdot\|\) denotes the 2-norm, the superscript \({}^{\top}\) represents the transpose of a matrix or vector, and \(f(z)\) is an analytic function in a given region and makes the corresponding matrix function well-defined \cite{higham2008}. The expression \eqref{1.1} arises from many real-world applications, such as option pricing based on the continuous-time Markov chain (CTMC) method \cite{Cui2021}, Markov PERT network \cite{Burgelman19}, computational electrodynamics \cite{Lambers10}, complex networks \cite{Benzi13,Fenu2017}, tuning scientific and probabilistic machine learning models \cite{NEURIPS2024}, etc. One is typically interested only in an approximation of Eq. \eqref{1.1}; thus, developing efficient algorithms to evaluate the bilinear form associated with a matrix function has become of paramount importance.  

The computation of bilinear forms of matrix functions has been extensively studied over the past several decades. Dahlquist et al. \cite{Dahlquist1972} first proposed the study of error bounding for linear systems in 1972. In 1994, Golub and Strako\v{s}  \cite{Golub1994} further investigated numerical algorithms for the quadratic form $\mathbf{u}^\top f(A)\mathbf{u}$ and extended them to the computation of the bilinear form \eqref{1.1}. Later, the work \cite{Bai1996} connected the numerical solution of the bilinear form of matrix functions with orthogonal polynomials, matrix traces, and quadrature formulas, proposed an integral solution process for the expression \eqref{1.1}, and provided upper and lower bounds for this problem under the somewhat restrictive condition that the high-order derivatives of matrix functions are computable. 



To evaluate the matrix function bilinear form \eqref{1.1} for a symmetric matrix $A$, a direct approach is to utilize the spectral decomposition $A = Q\Lambda Q^\top$. However, this direct method is restricted to symmetric settings and becomes computationally prohibitive for large and sparse matrices. To address these limitations, several efficient alternative approximations have been developed. These primarily include integral-based quadrature methods that reformulate the bilinear form into Riemann-Stieltjes integrals \cite{ Golub1994-1, Golub2010-1}, as well as Krylov subspace projection methods that transform the high-dimensional evaluation into a low-dimensional approximation problem \cite{Bai1996,Golub1997}.

As noted in the real applications mentioned previously, the involved matrix $A$ is typically nonsymmetric, for which the bi-Lanczos method is suitably employed to construct Krylov subspace based approximations. Calvetti et al.  \cite{calvetti1999applications} followed this idea to introduce biorthogonal polynomials to handle such nonsymmetric matrices, this bi-Lanczos framework has seen numerous advancements. Subsequent developments include Reichel et al.'s generalized averaged Gauss quadrature \cite{Reichel2015}, Schweitzer's two-sided extended Krylov method \cite{Schweitzer17}, and Pozza et al.'s theoretical unification of the Lanczos algorithm with complex Gauss quadrature \cite{Pozza2018Lanczos}. More recently, Alahmadi et al. \cite{alahmadi21} introduced Gauss-Laurent quadrature rules, and Altheimer \cite{altheimer2021} enhanced computational stability by treating (\ref{1.1}) as perturbed quadratic forms. However, the bi-Lanczos process requires matrix-vector products with both $A$ and $A^\top$ and is prone to breakdown. While look-ahead strategies \cite{fr2,Guo2003} mitigate serious breakdowns, they sacrifice short recurrences and increase computational costs. Alternatively, the Arnoldi process is widely used to approximate these bilinear forms \cite{Calvetti2005}, with Jia and Sun \cite{jia2020} providing a reliable a posteriori error estimate. However, its orthogonalization cost grows linearly with iterations, heavily consuming storage and limiting scalability.

Motivated by the above observations, we seek a more cost-effective alternative and intend to investigate the use of the Induced Dimension Reduction  (IDR($s$)) method for problem \eqref{1.1}. The IDR($s$) method is an iterative algorithm widely used for solving the system of linear equations. Since its introduction, researchers have developed various variants, such as IDR($s$)-biortho \cite{VanGijzen2011}, QMR-IDR($s$) \cite{Du2011a,Gijzen2014} and IDR($s$)-Ritz \cite{Simoncini10}. In practical applications, Sangers and van Gijzen \cite{SANGERS2015} successfully applied IDR($s$) to detect link-based spam in web graphs by solving large-scale singular linear systems derived from the Google matrix. Furthermore, IDR($s$) has been extended to eigenvalue computation using a generalized Hessenberg decomposition \cite{Zemke10,Gutknecht2013}. Additionally, these developments encourage numerical methods for other matrix computational problems \cite{rendel2013}. However, such extensions cannot directly yield high-precision eigenvectors for the standard eigenvalue problem. To overcome this limitation, Astudillo and van Gijzen \cite{Astudillo2019} proposed an improved variant based on a standard Hessenberg decomposition, which enhanced the accuracy of eigenvector approximations while maintaining computational efficiency, and further extended the approach to solving quadratic eigenvalue problems \cite{Astudillo17}. Motivated by recent advances in standard Hessenberg decomposition, this paper focuses on adapting the IDR($s$) method to compute the bilinear form involving matrix functions \eqref{1.1}. Additionally, we discuss the reliable stopping criteria for the proposed algorithm.

The rest of this paper is organized as follows. Section \ref{sec:main} reviews the basic principles of the IDR($s$) method. Section \ref{sec:err} introduces the IDR($s$)-based approximation for problem \eqref{1.1}. We then derive an explicit error expansion to provide a computable a posteriori error indicator for analytic matrix functions. Section \ref{sec:example} presents numerical experiments that assess the reliability of the proposed error indicator for several commonly used matrix functions, including the exponential, trigonometric, and square-root functions. Finally, concluding remarks are given in Section \ref{sec:conclusions}.
\section{The IDR($s$) method}
\label{sec:main}
This section reviews the derivation for constructing standard Hessenberg decompositions using the IDR($s$) method. Indeed, the IDR($s$) method is based on the so-called {\em IDR theorem}. 
Here, the theoretical underpinning of this method is encapsulated in the following theorem:
\begin{theorem}{\rm(IDR theorem \cite{Sonneveld2009})}
Let $A \in \mathbb{R}^{n\times n}$, $P = [\mathbf{p}_1, \mathbf{p}_2, \dots, \mathbf{p}_s] \in \mathbb{R}^{n\times s}$, and let $\{\mu_j\}$ be a sequence of scalars in $\mathbb{R}$. Defining the initial subspace as $\mathcal{G}_0 \equiv \mathbb{R}^n$, the subsequent subspaces are generated via the recurrence:
    \begin{equation*}
        \mathcal{G}_{j+1} = (A - \mu_{j+1} I)\left(\mathcal{G}_j\cap P^{\bot}\right),\quad j=0,1,\ldots,
    \end{equation*}
    where $P^{\bot}$ denotes the orthogonal complement of $P$. If $P^{\bot}$ does not contain any eigenvectors of $A$, then for all $j = 0, 1, \ldots$, the following hold:
     \begin{itemize}
        \item[1)] $\mathcal{G}_{j+1}\subset\mathcal{G}_j$,
        \item[2)] ${\rm dim}(\mathcal{G}_{j+1}) < {\rm dim}(\mathcal{G}_j)$, unless $\mathcal{G}_j = \{{\mathbf{0}}\}$.
     \end{itemize}
\end{theorem}

In early implementations of the IDR($s$) algorithm, the implicit nature of the recurrence precluded the explicit construction of a Hessenberg decomposition \cite{Sonneveld2009,VanGijzen2011}. Subsequently, Gutknecht and Zemke \cite{Gutknecht2013} (alongside related work in \cite{Gijzen2014}) derived a generalized Hessenberg decomposition from the IDR($s$) framework; we recall its derivation here. To explicitly construct the basis matrix $V_m$, an intermediate vector $\tilde{\mathbf{v}}_{i+1}$ residing in $\mathcal{G}_j$ can be expressed as:
\begin{equation*}
    \tilde{\mathbf{v}}_{i+1} = \left(A-\mu_{j+1} I\right)\left(\mathbf{v}_i-\sum_{\ell=1}^s c_{\ell} \mathbf{v}_{i-\ell}\right),
\end{equation*}
Normalize $\tilde{\mathbf{v}}_{i+1}$ to have unit length:
\begin{equation*}
    \tilde{\mathbf{v}}_{i+1} = \beta_{i+1} \mathbf{v}_{i+1},
\end{equation*}
where the $s+1$ preceding vectors $\mathbf{v}_{i-s}, \mathbf{v}_{i-s+1}, \ldots, \mathbf{v}_i$ belong to $\mathcal{G}_{j-1}$, $\mu_j \in \mathbb{R}$, and $j = \lfloor\frac{i}{s+1}\rfloor$. The coefficients $c_{\ell}$ are uniquely determined by solving the local $s \times s$ linear system:
\begin{equation}
    \left(P^\top\left[\mathbf{v}_{i-s}, \mathbf{v}_{i-s+1}, \ldots, \mathbf{v}_{i-1}\right]\right)\mathbf{c} = P^\top \mathbf{v}_i,
\label{eq:idr_system}
\end{equation}
where $\mathbf{c} = [c_1, c_2,\ldots,c_s]^\top$. Expanding the relation using Eq. \eqref{eq:idr_system} yields:
\begin{equation*}
    A \mathbf{v}_i = \beta_{i+1}\mathbf{v}_{i+1} + \mu_{j+1}\mathbf{v}_i - \mu_{j+1} \sum_{\ell=1}^s c_{\ell} \mathbf{v}_{i-\ell} + \sum_{\ell=1}^s c_{\ell} A\mathbf{v}_{i-\ell},
\end{equation*}
which can be equivalently rewritten to isolate the operations involving $A$:
\begin{equation*}
    A \mathbf{v}_i - \sum_{\ell=1}^s c_{\ell} A\mathbf{v}_{i-\ell} = \beta_{i+1}\mathbf{v}_{i+1} + \mu_{j+1}\mathbf{v}_i - \mu_{j+1} \sum_{\ell=1}^s c_{\ell} \mathbf{v}_{i-\ell}.
\end{equation*}

In the literature \cite{Gutknecht2013,Astudillo2016}, this sequence of equalities is encapsulated within the generalized Hessenberg decomposition:
\begin{equation}
    A V_m U_m = V_m \hat{H}_m + \beta_{m+1}\mathbf{v}_{m+1}\mathbf{e}_m^\top,
\label{eq:generalized_hessenberg}
\end{equation}
where $U_m$ is an upper triangular matrix and $\hat{H}_m$ is an upper Hessenberg matrix, forming the Sonneveld pencil $(\hat{H}_m, U_m)$. However, using the decomposition \eqref{eq:generalized_hessenberg} for matrix function approximation is computationally inefficient due to extra matrix manipulations. Fortunately, an equivalent \textit{standard} Hessenberg decomposition follows from the IDR($s$) recurrence \cite{Astudillo2016}. Exploiting the inherent structure of the generated Krylov basis, $A \mathbf{v}_{i-\ell}$ can be systematically expressed as a linear combination of the preceding basis vectors, yielding $A \mathbf{v}_{i-\ell} = V_{i-\ell+1} \mathbf{h}_{i-\ell}$. Substituting this relation into the expanded IDR sequence directly produces the standard factorization:
\begin{equation}
\begin{split}
    AV_m & = V_{m+1}\bar{H}_m \\
    & = V_m H_m + h_{m+1,m}\mathbf{v}_{m+1}\mathbf{e}_m^\top.
\end{split}
\label{2.6}
\end{equation}
Crucially, the columns of the standard Hessenberg matrix $H_m$ can be updated recursively, obviating the need to explicitly assemble the upper triangular matrix $U_m$, noting that $h_{m+1,m}=\beta_{m+1}$. Furthermore, theoretical analysis shows that $H_m$ satisfies a strict algebraic identity with its generalized counterpart \cite{Astudillo2016}:
\begin{equation*}
    H_m = \hat{H}_m U_m^{-1}.
\end{equation*}
This standard form is readily extensible to a broader range of matrix computation challenges. For completeness, Algorithm \ref{Alg1} outlines the IDR($s$) process \cite{Astudillo2016}.
\begin{algorithm}[!htp]
\caption{IDR($s$) process}
\begin{algorithmic}[1]
\label{Alg1}
\STATE \textbf{Input:} System matrix $A \in \mathbb{R}^{n \times n}$, target iteration $m$, parameter $s \in \mathbb{N}$, shadow matrix $P \in \mathbb{R}^{n \times s}$, initial basis matrix $V_s \in \mathbb{R}^{n \times (s+1)}$, and an upper Hessenberg matrix $\bar{H}_s$ satisfying $A V_s = V_{s+1}\bar{H}_s$.
\STATE \textbf{Output:} Updated basis matrix $V_{m+1} \in \mathbb{R}^{n \times (m+1)}$ and upper Hessenberg matrix $\bar{H}_m \in \mathbb{R}^{(m+1) \times m}$ such that $A V_m = V_{m+1} \bar{H}_m$. 
\FOR{\(i=s+1,\ldots,m\)}
    \IF{$i$ is a multiple of \(s+1\)}
        \STATE Select parameter \(\mu_j\) for subspace \(\mathcal{G}_j\)
        \STATE $k = 0$ \hfill \COMMENT{Reset the vector counter within the current $\mathcal{G}_j$}
    \ENDIF
    \STATE Solve the \(s\times s\) linear system:
    \begin{equation*}
    \left(P^{\top}[{\mathbf{v}}_{i-s},\mathbf{ v}_{i-s+1},\ldots,\mathbf{v}_{i-1}]\right)\mathbf{c} = P^{\top}\mathbf{v}_i
    \end{equation*}
    \STATE \(\mathbf{v} = \mathbf{v}_i - \sum^{s}_{\ell=1}c_{\ell}\mathbf{v}_{i-\ell}\) \hfill \COMMENT{where \(\mathbf{v}\in \mathcal{G}_{j-1}\cap P^{\bot}\)}
    \STATE \(\mathbf{ v}_{i+1} = (A - \mu_j I)\mathbf{v}\) \hfill \COMMENT{Generate new vector in \(\mathcal{G}_j\)}
    \STATE Generate the \((i+1)\)-th column of $H$ according to \cite[Eq. 9]{Astudillo2016}
    \FOR{\(\ell=1\to k\)}
        \STATE \(\beta_{i-\ell} = \mathbf{ v}^{\top}_i\mathbf{ v}_{i-\ell}\)
    \ENDFOR
    \STATE \(\beta_i = \|\mathbf{v}_{i+1} - \sum^{k}_{\ell=1}\beta_{i-\ell}\mathbf{ v}_{i-\ell}\|_2\)
    \STATE \(\mathbf{ v}_{i+1} = \left(\mathbf{ v}_{i+1} - \sum^{k}_{\ell=1}\beta_{i-\ell}\mathbf{ v}_{i-\ell}\right)/\beta_i\)
    \FOR{\(\ell = 1 \to k\)}
        \STATE \(h_{i-\ell,i} = h_{i-\ell,i} + \beta_{i-\ell}\)
    \ENDFOR
    \STATE \(h_{i+1,i} = \beta_i\)
    \STATE \(V_{i+1} = [\mathbf{ v}_1, \mathbf{ v}_2,\ldots,\mathbf{v}_i, \mathbf{ v}_{i+1}]\) \hfill \COMMENT{Update IDR decomposition}
    \STATE $k = k + 1$ \hfill \COMMENT{Increment the counter within $\mathcal{G}_j$}
\ENDFOR
\end{algorithmic}
\end{algorithm}

Notably, its initial matrix decomposition is computed via the $s$-step Arnoldi method \cite{arnoldi1951principle} using modified Gram-Schmidt orthogonalization. While an $m$-step Arnoldi process produces a comparable Hessenberg decomposition \eqref{2.6}, it typically incurs a higher computational cost than the $m$-step IDR($s$) approach \cite[Section 3.1]{Astudillo2016}. Throughout the remainder of this paper, we assume that Algorithm \ref{Alg1} runs for $m$ steps without breakdown. Under this standard algebraic assumption, the generated basis matrix $V_m \in \mathbb{R}^{n \times m}$ is guaranteed to have full column rank, which strictly implies that its Moore-Penrose pseudo-inverse $V_m^\dagger$ satisfies
\begin{equation*}
V_m^\dagger V_m = I_m.
\end{equation*}
\section{Error analysis of the proposed approximation}
\label{sec:err}
This section details the construction of the IDR($s$)-based approximation for \eqref{1.1} and provides a rigorous analysis of the corresponding error expansion. Specifically, we prove that the remainder term of this expansion vanishes under mild conditions, thereby validating the expansion itself. Based on this analysis, the leading computable term is used as a practical a posteriori error indicator for general analytic functions. We further formulate an efficient stopping criterion based on this indicator.
\subsection{The error expansion of \(E_m(f)\)}
Drawing upon the rationale of sketched Krylov subspace methods \cite{cortinovis2024}, which also generate a non-orthogonal basis, we project the original problem \eqref{1.1} onto the Krylov subspace via the relation \eqref{2.6}:
\begin{equation*}
	\begin{split}
\mathbf{u}^{\top}f(A)\mathbf{v} &\approx \beta \mathbf{u}^{\top}V_mV^{\dag}_m f(A)V_m \mathbf{e}_1\\
		& \approx \beta \mathbf{u}^{\top}V_m f(V^{\dag}_mAV_m)\mathbf{e}_1,\\
	\end{split}
\end{equation*}
where \(\beta = \|\mathbf{v}\|_2\) and \(\mathbf{v}_1 = \mathbf{v}/\beta\). For a non-orthogonal basis $V_m$, the exact Krylov compression of $A$
is given by $V_m^\dagger A V_m$ rather than the Hessenberg matrix $H_m$. However, the Hessenberg relation \eqref{2.6} implies that these two matrices differ only in their last column. Following the idea of cheaper
approximations for non-orthogonal Krylov bases, we therefore replace
$V_m^\dagger A V_m$ by $H_m$ and define the IDR($s$)-based approximation as
\[
F_m := \beta\mathbf{u}^\top V_m f(H_m) \mathbf{e}_1.
\]



Motivated by the theoretical framework of \cite{Saad1992analysis}, we define the corresponding error $E_m$:
\begin{equation*}
    E_m(f) = \mathbf{u}^{\top} f(A) \mathbf{v} - F_m. 
\end{equation*}
In this paper, let $S$ be a closed convex set containing the numerical ranges $\mathcal{F}(A)$ and $\mathcal{F}(H_m)$ of matrices $A$ and $H_m$. We assume that the function $f(z)$ is analytic on an open neighborhood $\mathcal{U}$ containing $S$~\cite{jia2020}. The sequence of points \(\{t_m\}_{m = 0}^{\infty}\) comes from $S$, and among them, the same points are arranged consecutively, that is, \(t_i = t_{i + 1} = \cdots = t_j\). The functions $\{\phi_j(t)\}$ are a sequence of analytic functions satisfying the following conditions:
\begin{equation}
\begin{cases}
		\phi_0(t)=f(t),\\
		\phi_{j+1}(t)=\frac{\phi_j(t)-\phi_j\left(t_j\right)}{t-t_j},& j \geq 0.
\end{cases}
\label{3.2}
\end{equation}

We now turn to the following asymptotic analysis of the approximation error expansion, the main result of this paper.
\begin{theorem} 
\label{Theorem 1}
Under the aforementioned analyticity assumption about $f(z)$, if either of the following conditions holds:
\begin{enumerate}
    \item[(i)] $f(z)$ satisfies the derivative-growth condition:
\begin{equation} \label{eq:growth_cond}
\max_{z \in S} \left| f^{(j)}(z) \right| \le M j! \rho^{-j}, \quad \forall j \ge 0,
\end{equation}
for some positive constants $M$ and $\rho > C_2$, or
    \item[(ii)] $f(z)$ represents a chosen single-valued analytic branch of a general analytic function, and the distance $r$ from $S$ to the boundary of its domain of analyticity (including any singularities and branch cuts) satisfies $r > C_2$ (where $C_2 = \|A\|_2 + T$ with $T = \max_{z \in S} |z|$),
\end{enumerate}
then for any given sequence of points $\{t_m\}_{m=0}^\infty$ in $S$, the error of the IDR($s$)-based approximation converges and can be expanded as:
\begin{equation}
\begin{split}
E_m(f) & = \mathbf{u}^\top f(A)\mathbf{v} - \beta \mathbf{u}^\top V_m f(H_m)\mathbf{e}_1 \\
& = \beta h_{m+1,m}\sum_{j=1}^\infty \mathbf{e}_m^\top \phi_j(H_m)\mathbf{e}_1 \mathbf{u}^\top p_{j-1}(A)\mathbf{v}_{m+1},
\end{split}
\label{3.3}
\end{equation}
where $p_0(t) = 1$, $p_j(t) = (t-t_0)(t-t_1)\cdots(t-t_{j-1})$ for $j \ge 1$, $t_i \in S$ for $i \ge 0$.
\end{theorem}
\begin{proof}
First, define the IDR($s$)-based approximation error of \(\phi_i(A) \mathbf{v}\) as follows:
\begin{equation}
		\mathbf{s}_m^i=\phi_i(A) \mathbf{v}-\beta V_m \phi_i\left(H_m\right) \mathbf{e}_1. \notag
	\end{equation}
    According to the recursive relation \eqref{3.2}, we have \(f(t) = f\left(t_0\right) + \left(t - t_0\right) \phi_1(t)\). Then, by applying the Hessenberg decomposition \eqref{2.6}, we obtain
    	\begin{align}
		\mathbf{u}^\top f(A) \mathbf{v}& =  f\left(t_0\right) \mathbf{u}^\top \mathbf{v
        }+\mathbf{u}^\top\left(A-t_0 I\right) \phi_1(A) \mathbf{v} \notag\\
		& = f\left(t_0\right) \mathbf{u}^\top\mathbf{v} +\mathbf{u}^\top\left(A-t_0 I\right)\left(\beta V_m \phi_1\left(H_m\right) \mathbf{e}_1+\mathbf{s}_m^1\right) \notag\\
		& = f\left(t_0\right) \mathbf{u}^\top \mathbf{v}+\beta \mathbf{u}^\top\left(A-t_0I\right) V_m\phi_1\left(H_m\right) \mathbf{e}_1+\mathbf{u}^\top\left(A-t_0 I\right) \mathbf{s}_m^1\notag\\
		& = f\left(t_0\right) \mathbf{u}^\top \mathbf{v}+\beta \mathbf{u}^\top(V_m(H_m-t_0I)+h_{m+1,m}\mathbf{v}_{m+1}\mathbf{e}^{\top}_m)\phi_1(H_m) \mathbf{e}_1+\mathbf{u}^\top\left(A-t_0 I\right) \mathbf{s}_m^1\notag \\
		& = \beta \mathbf{u}^\top V_m\left(f\left(t_0\right) \mathbf{e}_1+\left(H_m-t_0 I\right) \phi_1\left(H_m\right) \mathbf{e}_1\right)+\beta h_{m+1,m}\mathbf{u}^{\top}\mathbf{v}_{m+1}\mathbf{e}^{\top}_m\phi_1(H_m)\mathbf{e}_1\notag\\
        &\quad + \mathbf{u}^\top\left(A-t_0 I\right) \mathbf{s}_m^1 \notag\\
		&= \beta \mathbf{u}^\top V_m f(H_m)\mathbf{e}_1+\beta  h_{m+1,m} \mathbf{u}^{\top}\mathbf{v}_{m+1}\mathbf{e}^{\top}_m\phi_1\left(H_m\right) \mathbf{e}_1+\mathbf{u}^\top\left(A-t_0 I\right) \mathbf{s}_m^1.\notag
	\end{align}
%
By applying the analogous algebraic manipulations to the function $\phi_1(A)\mathbf{v}$ with the recurrence relation \eqref{3.2}, we can obtain:
\begin{equation*}
\phi_1(A)\mathbf{v} = \beta  V_m \phi_1(H_m)\mathbf{e}_1 + \beta  h_{m+1,m} \mathbf{v}_{m+1}\mathbf{e}_m^\top \phi_2(H_m)\mathbf{e}_1 + (A-t_1 I)\mathbf{s}_m^2.
\end{equation*}
    
    Since \(\mathbf{s}_m^i = \phi_i(A)\mathbf{v} - \beta V_m\phi_i\left(H_m\right)\mathbf{e}_1\), substituting \(i = 1\) and the above results gives:
    \begin{equation}
    \begin{split}
		\mathbf{s}_m^1 & =\phi_1(A) \mathbf{v}-\beta V_m \phi_1\left(H_m\right) \mathbf{e}_1 \\
        & =\beta  h_{m+1,m}\mathbf{v}_{m+1}\mathbf{e}^{\top}_m\phi _2(H_m)\mathbf{e}_1+\left(A-t_1 I\right) \mathbf{s}_m^2\notag.
    \end{split}
	\end{equation}
    Similarly, for \(\mathbf{s}_m^2\), \(\mathbf{s}_m^3\), \(\cdots\), the general term has the following form:
    	\begin{equation}
		\mathbf{s}_m^{i-1}=\beta  h_{m+1,m}\mathbf{v}_{m+1}\mathbf{e}^{\top}_m\phi_i(H_m) \mathbf{e}_1+\left(A-t_{i-1} I\right) \mathbf{s}_m^i, \quad i=2,3, \cdots, \infty\notag.
	\end{equation}
    Substituting the above general formula into the available equations yields:
    	\begin{align}
		E_m(f) & = \mathbf{u}^{\top}f(A)\mathbf{v}-\beta \mathbf{u}^{\top}V_mf(H_m)\mathbf{e}_1 \notag\\
		& =  \beta  h_{m+1,m}\mathbf{u}^{\top}\mathbf{v}_{m+1}\mathbf{e}^{\top}_m\phi_1(H_m)
		\mathbf{e}_1+\mathbf{u}^{\top}(A-t_0I)\mathbf{s}_m^1 \notag\ \\
		& = \beta  h_{m+1,m}\sum_{j=1}^{i}\mathbf{e}^{\top}_m\phi_j(H_m)\mathbf{e}_1\mathbf{u}^{\top}p_{j-1}(A){\bf v}_{m+1}+\mathbf{u}^{\top}p_i(A)\mathbf{s}_m^i\notag.
		\end{align}

The convergence rate analysis of $|\mathbf{u}^\top p_i(A)\mathbf{s}^i_m|$ is given below. We show that the remainder term converges geometrically to zero under the assumptions of Theorem \ref{Theorem 1}. From the definition of the function sequence $\{\phi_i(t)\}$, $\phi_{j+1}(t)$
can be expressed by a divided difference as
\begin{equation}
	\phi_{j+1}(t)=f[t,t_0,t_1,\ldots,t_j]. \notag
\end{equation}
Since \( f(z) \) is analytic on the closed convex set \( S \) containing the nodes, by using the Hermite-Genocchi formula for the divided difference, we have:
\[
f[t, t_0, \dots, t_{i-1}] = \int_{\Delta^i} f^{(i)}(\tau_0 t + \tau_1 t_0 + \dots + \tau_i t_{i-1}) \, d\tau.
\]

Since the matrix $A$ is typically non-symmetric, its eigenvalues and numerical range generally lie in the complex plane $\mathbb{C}$. To establish rigorous bounds in the complex domain,  we utilize the growth condition \eqref{eq:growth_cond}, then the standard bound for the $j$-th order divided difference in Eq. \eqref{3.3} can be rewritten as:
\begin{equation} \label{eq:divided_diff_new}
|f[z_0, \dots, z_j]| \le \frac{\max_{z\in S} |f^{(j)}(z)|}{j!} \le M \rho^{-j}.
\end{equation}
By virtue of Eq. \eqref{eq:divided_diff_new}, we can immediately obtain the bounds for both the function $\phi_i(t)$ and its $k$-th derivative $\phi_i^{(k)}(z)$.
By the properties of divided differences, the $k$-th derivative of $\phi_i(t)$ can be represented as a divided difference of order $k+i$:
\begin{equation*}
\phi_i^{(k)}(z) = k! f[\underbrace{z, \dots, z}_{k+1}, t_0, \dots, t_{i-1}].
\end{equation*}
Therefore, we obtain the geometrically decaying bounds for $\phi_i(t)$ and its derivatives:
\begin{equation} \label{eq:phi_deriv_new}
\sup_{z \in S} \left| \phi_i^{(k)}(z) \right| \le k! \frac{\max_{\zeta \in S} \left| f^{(k+i)}(\zeta) \right|}{(k+i)!} \le k! \frac{M (k+i)! \rho^{-(k+i)}}{(k+i)!} = M k! \rho^{-(k+i)}.
\end{equation}
Eq. \eqref{eq:phi_deriv_new} justifies the estimations for $\phi_i(t)$ and its derivatives in the complex setting without requiring the nodes to be real.
Based on these estimations, we can complete the convergence analysis of the remainder term $\left| \mathbf{u}^\top p_i(A) \mathbf{s}_m^i \right|$ in Theorem \ref{Theorem 1}. Appendix A proves that under the conditions of Theorem \ref{Theorem 1}, the norm of the error vector satisfies
\begin{equation*}
\|\mathbf{s}_m^i\| \le C_\rho \rho^{-i},
\end{equation*}
and the remainder term can be strictly bounded by
\begin{equation*}
\left| \mathbf{u}^\top p_i(A) \mathbf{s}_m^i \right| \le \hat{C}_\rho \left( \frac{C_2}{\rho} \right)^i,
\end{equation*}
where $\rho > C_2$. To keep the readability and flow of the main text, the detailed, rigorous mathematical proof (covering both entire and general analytic functions) is deferred to Appendix \ref{app:convergence_proof}. This completes the sketch of the proof.
\end{proof}

\begin{remark} \label{rem:analytic_convergence}
For general analytic functions with singularities (such as $z^{-1}$ and $z^{1/2}$), the divided differences and their derivatives exhibit geometric decay rather than superlinear decay. Under the assumption that the domain of analyticity of $f(z)$ is sufficiently large, the remainder term $\mathbf{u}^\top p_i(A)\mathbf{s}_m^i$ is guaranteed to converge geometrically to zero. This mathematically ensures the convergence of the  resulting error expansion \eqref{3.3} for general analytic functions without relying on Stirling's approximation \cite{jia2020}. The detailed, rigorous mathematical proof is presented in Appendix~\ref{app:general_analytic_proof}.
\end{remark}
  
 Note that if the nodes $\{t_i\}_{i=0}^{n-1}$ are the $n$ eigenvalues of the matrix $A$, 
 then the resulting error expansion \eqref{3.3} reduces to a finite sum, and
 \begin{equation*}
\begin{split}
E_m(f) & = \mathbf{u}^{\top}f(A)\mathbf{v}-\beta \mathbf{u}^\top V_mf(H_m)\mathbf{e}_1 \\
& = \beta h_{m+1,m} \sum_{j=1}^{n}\mathbf{e}^{\top}_m\phi_j(H_m)\mathbf{e}_1\mathbf{u}^\top p_{j-1}(A)\mathbf{v}_{m+1}.
\end{split}
\end{equation*}
	 \begin{remark}
     \vspace{-1em}
     The assumption in Theorem \ref{Theorem 1} that $S$ contains both $\mathcal{F}(A)$ and $\mathcal{F}(H_m)$ warrants clarification. Unlike the standard Arnoldi process, where the orthogonal projection guarantees $\mathcal{F}(H_m) \subseteq \mathcal{F}(A)$, the IDR($s$) method employs an oblique projection via a non-orthogonal basis $V_m$. Quantitatively, the norm of the projected matrix satisfies $\|H_m\|_2 \le \kappa_2(V_m) \|A\|_2 + |h_{m+1, m}| \|V_m^\dagger\|_2$, where $\kappa_2(V_m) = \|V_m\|_2 \|V_m^\dagger\|_2$ is the condition number. If the basis $V_m$ becomes highly ill-conditioned (i.e., $\kappa_2(V_m) \gg 1$), $\mathcal{F}(H_m)$ can expand significantly, generating spurious Ritz values could theoretically deviate towards the singularities of $f(z)$. Thus, explicitly assuming that $S$ covers $\mathcal{F}(H_m)$ is imperative to validate the contour integral bounds in our error analysis. Nevertheless, empirical observations indicate that, in the absence of severe breakdown, $\mathcal{F}(H_m)$ remains stably localized around the spectrum of $A$. Alternatively, the expansion of the numerical range could be mitigated by adapting the similarity-restoring correction recently proposed for randomized Krylov subspace methods \cite{grigori2026}. By enforcing orthogonality on the last basis vector, this technique strictly controls the spectral properties of the projected matrix. We will explore this in greater depth in future research.
  \end{remark}
%

\begin{remark} 
\label{rem:tail_sum_rigorous}
For the matrix exponential $f(z) = e^z$, which is an entire function, for every fixed $\rho > C_2$, there exists a constant $M_\rho > 0$ such that the growth condition \eqref{eq:growth_cond} holds. To ensure a rapid decay rate of the remainder term, for any fixed radius $\rho > C_2$ for which the derivative-growth condition holds, the tail admits a geometric bound.
By the relation $R_{i-1} = a_i + R_i$, and since $R_i \to 0$ as $i \to \infty$ by the convergence analysis in Appendix \ref{eq:complex_schur}, the tail of the error series from $j = P$ to $\infty$ is mathematically equivalent to the remainder $R_{P-1}$:
\begin{equation} \label{eq:tail_sum_eq_remainder}
\beta h_{m+1,m} \sum_{j=P}^\infty \mathbf{e}_m^\top \phi_j(H_m)\mathbf{e}_1 \mathbf{u}^\top p_{j-1}(A)\mathbf{v}_{m+1} = R_{P-1},
\end{equation}
where $a_j = \beta h_{m+1, m} \mathbf{e}_m^\top \phi_j(H_m) \mathbf{e}_1 \mathbf{u}^\top p_{j-1}(A) \mathbf{v}_{m+1}$ represents the $j$-th term of the error series. Combining Eq. \eqref{eq:tail_sum_eq_remainder} with the rigorous bound derived in Eq. \eqref{eqA3} of Appendix \ref{eq:complex_schur}, and strictly preserving the constant factor, we obtain:
\begin{align*} 
\left| \beta h_{m+1,m} \sum_{j=P}^\infty \mathbf{e}_m^\top \phi_j(H_m)\mathbf{e}_1 \mathbf{u}^\top p_{j-1}(A)\mathbf{v}_{m+1} \right| &= \left| R_{P-1} \right| \notag \\
&\le C_4 \left( \frac{C_2}{\rho} \right)^{P-1},
\end{align*}
where the constant factor $C_4$ is strictly preserved. To ensure that the tail sum is bounded by a given target tolerance $\varepsilon > 0$, it is sufficient to choose the truncation index $P$ such that $C_4 (C_2 / \rho)^{P-1} \le \varepsilon$, which yields:
\begin{equation*} 
P \ge 1 + \frac{\log(C_4/\varepsilon)}{\log(\rho/C_2)}.
\end{equation*}
This geometric tail bound completely obviates the need for Stirling's approximation and the arbitrary threshold $P = \max\{1, \lceil e^2 C_2 \rceil\}$.

Numerical experiments suggest that the first few terms often dominate the error expansion. For the matrix exponential $e^{-hA}$, Proposition~\ref{prop:asymptotic_error} rigorously shows that, under the stated non-cancellation condition, the first term is asymptotically equivalent to the total error as $h \to 0$. For general analytic matrix functions, the leading term is therefore employed as a practical heuristic indicator, whose effectiveness is further examined numerically.
\end{remark}

\subsection{A computable error indicator and stopping criterion}
The error expansion in Theorem \ref{Theorem 1} provides a natural way to construct a computable error indicator for the IDR($s$)-based approximation. For any positive integer $q$, we define the $j$-th term of the error series as:
\begin{equation*} 
a_{j,m} = \beta h_{m+1,m} \mathbf{e}_m^\top \phi_j(H_m)\mathbf{e}_1 \mathbf{u}^\top p_{j-1}(A)\mathbf{v}_{m+1}.
\end{equation*}
Then, the true error $E_m(f)$ in Eq. \eqref{3.3} can be written as the finite-sum expansion:
\begin{equation*} 
E_m(f) = \sum_{j=1}^q a_{j,m} + R_{q,m}, \quad \text{with} \quad R_{q,m} = \mathbf{u}^\top p_q(A) \mathbf{s}_m^q.
\end{equation*}
Under the assumptions of Theorem \ref{Theorem 1} (and the convergence analysis in Appendix \ref{app:convergence_proof}), the remainder term satisfies $\lim_{q \to \infty} R_{q,m} = 0$. Hence, the first term:
\begin{equation*}
a_{1,m} = \beta h_{m+1,m} \mathbf{e}_m^\top \phi_1(H_m)\mathbf{e}_1 \mathbf{u}^\top \mathbf{v}_{m+1}
\end{equation*}
represents the leading computable contribution to the approximation error. We therefore define the absolute a posteriori error indicator $\eta_m$ as:
\begin{equation*} 
\eta_m = \beta |h_{m+1,m}| \left| \mathbf{e}_m^\top \phi_1(H_m)\mathbf{e}_1 \right| \left| \mathbf{u}^\top \mathbf{v}_{m+1} \right|.
\end{equation*}
Unlike a rigorous error upper bound, $\eta_m$ should be interpreted as an asymptotic error indicator. It is expected to approximate the absolute error $|E_m(f)|$ when the higher-order terms are sufficiently small and no severe cancellation occurs.

To rigorously justify this heuristic indicator for the ubiquitous matrix exponential $f_h(z) = e^{-hz}$ with step size $h \ll 1$, we specialize the analysis to the scaling parameter $h \to 0$. By the definition of the $\phi$-functions, we have the Taylor expansion $\phi_j(z) = \sum_{k=0}^\infty \frac{z^k}{(k+j)!}$ \cite{Saad1992}. Since $H_m \in \mathbb{R}^{m \times m}$ is an upper Hessenberg matrix, its powers satisfy the standard identity $\mathbf{e}_m^\top H_m^k \mathbf{e}_1 = 0$ for $k = 0, 1, \dots, m-2$. By the relation $\phi_j^{(h)}(z) = (-h)^j \varphi_j(-hz)$, the $j$-th term of the error series can be explicitly written as $$a_{j,m}(h) = \beta h_{m+1, m} (-h)^j \mathbf{e}_m^\top \varphi_j(-h H_m) \mathbf{e}_1 \mathbf{u}^\top A^{j-1} \mathbf{v}_{m+1}.$$ We now establish the following formal proposition regarding the asymptotic behavior of the matrix exponential error series:

\begin{proposition} \label{prop:asymptotic_error}
Let the $m$-step IDR($s$) decomposition of $A$ be fixed and assume that it proceeds without breakdown. For the matrix exponential $f_h(z) = e^{-hz}$, let $E_m(e^{-hA}) = a_{1,m}(h) + E_m^{(2)}(h)$. Then, as $h \to 0$,
\begin{equation} \label{eq:a1m_asymptotic}
a_{1,m}(h) = c_m h^m + \mathcal{O}_1(h^{m+1}),
\end{equation}
where
\begin{equation*} 
c_m = \frac{(-1)^m \beta h_{m+1, m}}{m!} \mathbf{e}_m^\top H^{m-1}_m \mathbf{e}_1 \mathbf{u}^\top \mathbf{v}_{m+1}.
\end{equation*}
Moreover, the higher-order remainder satisfies:
\begin{equation}\label{eq:rem_h}
E_m^{(2)}(h) = \mathcal{O}(h^{m+1}).
\end{equation}
If $\mathbf{u}^\top \mathbf{v}_{m+1} \neq 0$, then $c_m \neq 0$, and consequently, we obtain the strict error-estimator equivalence:
\begin{equation}\label{eq:strict_equivalence}
E_m(e^{-hA}) = a_{1,m}(h)\big(1 + \mathcal{O}(h)\big), \quad h \to 0.
\end{equation}
In particular, this implies the asymptotic exactness of the error indicator:
\begin{equation*}
\lim_{h \to 0} \frac{|E_m(e^{-hA})|}{|a_{1,m}(h)|} = 1.
\end{equation*}
\end{proposition}
\begin{proof}
For any fixed $j \ge 1$, as $h \to 0$, we expand $\mathbf{e}_m^\top \varphi_j(-hH_m)\mathbf{e}_1$ using the upper Hessenberg matrix property $\mathbf{e}_m^\top H_m^k \mathbf{e}_1 = 0$ for $k < m-1$:
\begin{align} \label{eq:phi_asymptotic_order}
\mathbf{e}_m^\top \varphi_j(-hH_m)\mathbf{e}_1 &= \sum_{k=m-1}^\infty \frac{(-h)^k}{(k+j)!} \mathbf{e}_m^\top H_m^k \mathbf{e}_1 \notag \\
&= \frac{(-h)^{m-1}}{(m+j-1)!} \mathbf{e}_m^\top H_m^{m-1} \mathbf{e}_1 + \mathcal{O}_j(h^m),
\end{align}
where the notation $\mathcal{O}_j$ indicates that the implied constant may depend on $j$. For $j=1$, substituting Eq. \eqref{eq:phi_asymptotic_order} into the definition of $a_{1,m}(h)$ yields Eq. \eqref{eq:a1m_asymptotic}.

To rigorously prove Eq. \eqref{eq:rem_h} without relying on any unproved uniform assumptions on the $\mathcal{O}_j$ constants, we establish a uniform bound on $a_{j,m}(h)$ for any $0 < h \le h_0$, where $h_0$ is a fixed positive step size. Using the series expansion, we have:
\begin{equation} \label{eq:phi_series_strict}
\left| \mathbf{e}_m^\top \varphi_j(-h H_m) \mathbf{e}_1 \right| \le \sum_{k=m-1}^\infty \frac{h^k \|H_m\|_2^k}{(k+j)!} = h^{m-1} \sum_{\ell=0}^\infty \frac{h^\ell \|H_m\|_2^{m-1+\ell}}{(m-1+\ell+j)!},
\end{equation}
where we have substituted $k = m-1+\ell$. Since the binomial coefficient satisfies $\binom{m-1+\ell+j}{\ell} \ge 1$, we have the standard factorial inequality:
\begin{equation} \label{eq:factorial_ineq}
\frac{1}{(m-1+\ell+j)!} \le \frac{1}{(m+j-1)! \ell!}.
\end{equation}
Substituting Eq. \eqref{eq:factorial_ineq} into Eq. \eqref{eq:phi_series_strict} and using $h \le h_0$, we obtain:
\begin{equation*}
\left| \mathbf{e}_m^\top \varphi_j(-h H_m) \mathbf{e}_1 \right| \le \frac{h^{m-1}}{(m+j-1)!} \sum_{\ell=0}^\infty \frac{h_0^\ell \|H_m\|_2^{m-1+\ell}}{\ell!} = \frac{h^{m-1}}{(m+j-1)!} C_{m,0}(h_0),
\end{equation*}
where $C_{m,0}(h_0) = \|H_m\|_2^{m-1} e^{h_0 \|H_m\|_2}$ is a constant independent of $h$ and $j$. Substituting this back into the term definition, for any $j \ge 1$ and $0 < h \le h_0$, we have the uniform bound:
\begin{align} \label{eq:ajm_uniform_bound}
\left| a_{j,m}(h) \right| &= \beta |h_{m+1, m}| h^j \left| \mathbf{e}_m^\top \varphi_j(-hH_m) \mathbf{e}_1 \right| \left| \mathbf{u}^\top A^{j-1} \mathbf{v}_{m+1} \right| \notag \\
&\le C_m(h_0) \frac{h^{m+j-1} \|A\|_2^{j-1}}{(m+j-1)!},
\end{align}
where $C_m(h_0) = \beta |h_{m+1, m}| C_{m,0}(h_0)$ is a constant independent of $h$ and $j$. 

Using the uniform bound \eqref{eq:ajm_uniform_bound}, we sum the tail of the error series from $j=2$ to $\infty$ (substituting $\ell = j-2$):
\begin{align} \label{eq:em2_strict_sum}
\left| E_m^{(2)}(e^{-hA}) \right| &\le \sum_{j=2}^\infty \left| a_{j,m}(h) \right| \le C_m(h_0) \sum_{j=2}^\infty \frac{h^{m+j-1} \|A\|_2^{j-1}}{(m+j-1)!} \notag \\
&= C_m(h_0) h^{m+1} \|A\|_2 \sum_{\ell=0}^\infty \frac{h^\ell \|A\|_2^\ell}{(m+1+\ell)!} \notag \\
&\le C_m(h_0) h^{m+1} \|A\|_2 \sum_{\ell=0}^\infty \frac{h_0^\ell \|A\|_2^\ell}{(m+1)!\ell!} = \mathcal{O}(h^{m+1}),
\end{align}
which rigorously proves Eq. (3.11). Under the non-cancellation condition $\mathbf{u}^\top \mathbf{v}_{m+1} \neq 0$, we have $c_m \neq 0$ from Eq. (3.10), and thus $|a_{1,m}(h)| \asymp h^m$. Since the total error satisfies $E_m(e^{-hA}) = a_{1,m}(h) + E_m^{(2)}(e^{-hA})$, this, combined with the remainder bound $E_m^{(2)}(e^{-hA}) = \mathcal{O}(h^{m+1})$ in Eq. \eqref{eq:em2_strict_sum}, immediately yields:
\begin{equation*}
\frac{E_m(e^{-hA})}{a_{1,m}(h)} = 1 + \frac{E_m^{(2)}(e^{-hA})}{a_{1,m}(h)} = 1 + \mathcal{O}(h), \quad h \to 0.
\end{equation*}
Taking the absolute value directly leads to the strict error-estimator equivalence \eqref{eq:strict_equivalence}, completing the proof.
\end{proof}

To prevent premature termination caused by accidental stagnation or temporary cancellations, we combine the error indicator $\eta_m$ with the change in two consecutive approximations. We adopt a standard dual-tolerance criterion widely used in robust numerical libraries. The iteration is terminated if both the relative/absolute error indicator and the difference between consecutive approximations satisfy:
\begin{align} 
\max\{\eta_m, \eta_{m-1}\} &\le \texttt{atol} + \texttt{rtol} |F_m|, \label{eq:stop_criterion_rtol}\\
|F_m - F_{m-1}| &\le \tau \left( \texttt{atol} + \texttt{rtol} |F_m| \right), \label{eq:stop_criterion_rtol2}
\end{align}
where $\texttt{atol}$ and $\texttt{rtol}$ are the user-specified absolute and relative tolerances, respectively, and $\tau \ge 1$ is a safety scaling factor.

\begin{remark}\label{rem:dual_tolerance}
In practice, to prevent premature termination caused by accidental stagnation, the proposed error estimator $\xi_1^I$ can be easily embedded into a standard dual-tolerance criterion \eqref{eq:stop_criterion_rtol}-\eqref{eq:stop_criterion_rtol2}. However, since severe stagnation was not observed in our standard test problems, For the standard comparison experiments in Examples 1-4, we use the simple relative criterion to ensure a fair comparison among the three methods. Example 5 then separately investigates the robustness of the dual criterion under a controlled temporary cancellation.

\end{remark}

At the end of this section, we present the complete IDR($s$)-based algorithm for evaluating the bilinear forms of matrix functions in Algorithm \ref{Alg_IDR}. It is worth noting that since the first $s$ steps of the IDR($s$) process are just the Arnoldi process, the relative error curves of the two methods for computing the matrix function bilinear form are identical during these initial $s$ steps.
\begin{algorithm}[!htb]
\caption{Computing the bilinear form \eqref{1.1}}
\label{Alg_IDR}
\begin{algorithmic}[1]
\STATE \textbf{Input:} Matrix $A \in \mathbb{R}^{n \times n}$, projection vectors $\boldsymbol{u}, \boldsymbol{v} \in \mathbb{R}^n$, analytic function $f$, IDR parameter $s$, tolerance \textbf{tol}, and maximum number of iterations \textbf{maxit}.
\STATE \textbf{Output:} The approximated value  $F_m \approx \boldsymbol{u}^\top f(A) \boldsymbol{v}$.
\STATE \(\beta = \|\boldsymbol{v}\|_2\), \(\boldsymbol{v}_1 = \boldsymbol{v}/\beta\)
\FOR{\(k = 1\) to maxit}
    \STATE Generate the Hessenberg decomposition $AV_m = V_m H_m +  h_{m+1,m}\mathbf{v}_{m+1}{\mathbf{e}}^{\top}_m$ \hfill \COMMENT{Call Algorithm \ref{Alg1}}
    \STATE Compute the approximation: $F_m =\beta \boldsymbol{u}^\top (V_m \cdot f(H_m) \cdot \boldsymbol{e}_1)$ 
    \STATE Verify the a posteriori error estimate using the leading term of the  resulting error expansions in Remark \ref{rem:tail_sum_rigorous} or Theorem \ref{Theorem 1}.
    \IF{\(\mathrm{mode}=\mathrm{single}\) and \(\mathrm{relres}<\mathrm{tol}\)}
    \STATE Stop
\ELSIF{\(\mathrm{mode}=\mathrm{dual}\), \(m>1\), and criteria \eqref{eq:stop_criterion_rtol}--\eqref{eq:stop_criterion_rtol2} are satisfied}
    \STATE Stop
\ENDIF
\ENDFOR
\end{algorithmic}
\end{algorithm}
\section{Numerical Experiments}
\label{sec:example}
All the numerical experiments presented in this study were performed on a Windows 10 (64 bit) computer equipped with 16.0GB of 12th Gen Intel(R) Core(TM) i7-12700 processor operating at 2.10 GHz. MATLAB R2025a served as the software platform for all computations. 

In this section, we evaluate and compare the performance of the approximations based on three different processes: Arnoldi, IDR($s$), and bi-Lanczos. Table \ref{table111} summarizes their underlying Hessenberg decompositions, approximations ($F_m^*$), and a posteriori error estimate ($\xi_1^*$). Here, $V_m^*$ is the generated basis matrix, $H_m^*$ (or $T_m^*$) is the projection matrix, $\mathbf{v}_{m+1}^*$ is the next basis vector, and $h_{m+1,m}$ or $\delta_{m+1}$ are the normalization coefficients. To compute the functions $\phi_i(H_m^*)$ required for error estimation, a shift $t_0$ must be specified within the numerical range of $H_m^*$. Since our experiments indicate that the estimator is insensitive to $t_0$, we set $t_0 = 0$ or $h_{1,1}$. Given the moderate size of $H_m^*$, we directly utilize MATLAB's built-in \texttt{funm} to evaluate $f(H_m^*)$ and $\phi_i(H_m^*)$. However, for the matrix square root, the built-in function \texttt{sqrtm} will be used.
\begin{table}[!htb]
\centering
\caption{Matrix decompositions, approximations, and a posteriori error estimate for bilinear forms}
\label{table111}
\begin{adjustbox}{width=\textwidth}
\begin{tabular}{l c c c}
\toprule
\textbf{Method} & \textbf{Hessenberg Decomposition} & \textbf{Approximation} & \textbf{a posteriori error estimate} \\
\midrule
Arnoldi
& $\displaystyle A V_m^{A} = V_m^{A} H_m^{A} + h_{m+1,m}^{A} \mathbf{v}_{m+1}^{A} \mathbf{e}_m^\top$
& $\displaystyle F_m^{A} = \beta \mathbf{u}^\top V_m^{A} f(H_m^{A}) \mathbf{e}_1$
& $\displaystyle \xi_1^{A} = \beta |h_{m+1,m}^{A}| \left| \mathbf{e}_m^\top  \phi_1(H_m^{A}) \mathbf{e}_1 \right| \left| \mathbf{u}^\top \mathbf{v}_{m+1}^{A} \right|$ \\
IDR($s$)
& $\displaystyle A V_m^{I} = V_m^{I} H_m^{I} + h_{m+1,m}^I\mathbf{v}_{m+1}^{I} \mathbf{e}_m^\top$
& $\displaystyle F_m^{I} = \beta\mathbf{u}^\top V_m^{I} f(H_m^{I}) \mathbf{e}_1$
& $\displaystyle \xi_1^{I} = \beta  h_{m+1,m}^I\left| \mathbf{e}_m^\top  \phi_1(H_m^{I}) \mathbf{e}_1 \right| \left| \mathbf{u}^\top \mathbf{v}_{m+1}^{I} \right|$ \\
Bi-Lanczos
& $\begin{aligned}
    A V_m^{L} &= V_m^{L} T_m^{L} + \delta_{m+1} \mathbf{v}_{m+1}^{L} \mathbf{e}_m^\top 
  \end{aligned}$
& $\displaystyle F_m^{L} = \beta \mathbf{u}^\top V_m^{L} f(T_m^{L}) \mathbf{e}_1$
& $\displaystyle \xi_1^{L} = \beta |\delta_{m+1}| \left| \mathbf{e}_m^\top \phi_1(T_m^{L}) \mathbf{e}_1 \right| \left| \mathbf{u}^\top \mathbf{v}_{m+1}^{L} \right|$ \\
\bottomrule
\end{tabular}
\end{adjustbox}
\end{table}

To standardize the numerical evaluation, we define the true relative error as:
\begin{equation}
    \xi_{\text{true}}^{\text{rel}} = \frac{|\mathbf{u}^\top f(A) \mathbf{v} - F_m^*|}{|\mathbf{u}^\top f(A) \mathbf{v}|}.
\end{equation}
Because the exact value $\mathbf{u}^\top f(A) \mathbf{v}$ is practically inaccessible, we replace the denominator with the computed approximation $|F_m^*|$. Combined with the absolute error estimator $\xi_1^*$ from Table \ref{table111}, this yields the estimated relative error:
\begin{equation*}
    \xi^{\text{rel}} = \frac{\xi_1^*}{|F_m^*|}. 
\end{equation*}
In our experiments, the random vectors $\mathbf{u}$ and $\mathbf{v}$ are normalized such that $\|\mathbf{u}\|_2 = \|\mathbf{v}\|_2 = 1$ (implying $\beta = 1$). For all methods, the stopping tolerance is fixed at $10^{-8}$, with the IDR($s$) parameter chosen as $s=6$.
\textbf{Example 1.}
The test matrix $G \in \mathbb{R}^{n \times n}$ is the well-known \texttt{grcar} matrix of order $n \in \{10000, 15000, 20000\}$, which is a real, highly non-symmetric matrix sourced from Higham's Matrix Computation Toolbox\footnote{\url{https://www.mathworks.com/help/matlab/ref/gallery.html}.}. To comprehensively evaluate the numerical stability and convergence behavior of the algorithms, the scaling parameter $h$ is systematically chosen as \(h = 0.2, 0.5, 1\).

\begin{figure}[!htp]
    \centering

    \begin{subfigure}[b]{0.31\textwidth}
        \centering
        \includegraphics[width=\textwidth]{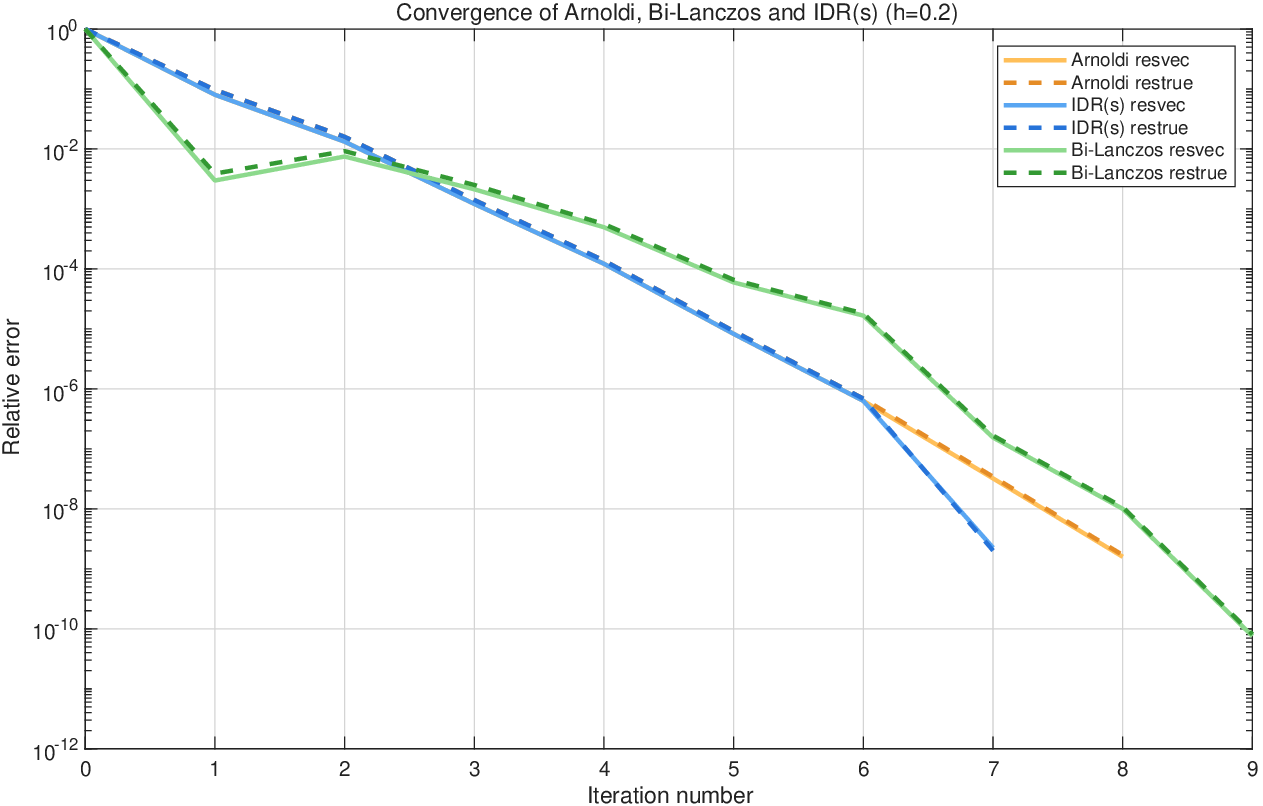}
        \caption*{\tt $h=0.2$}
    \end{subfigure}\hfill
    \begin{subfigure}[b]{0.31\textwidth}
        \centering
        \includegraphics[width=\textwidth]{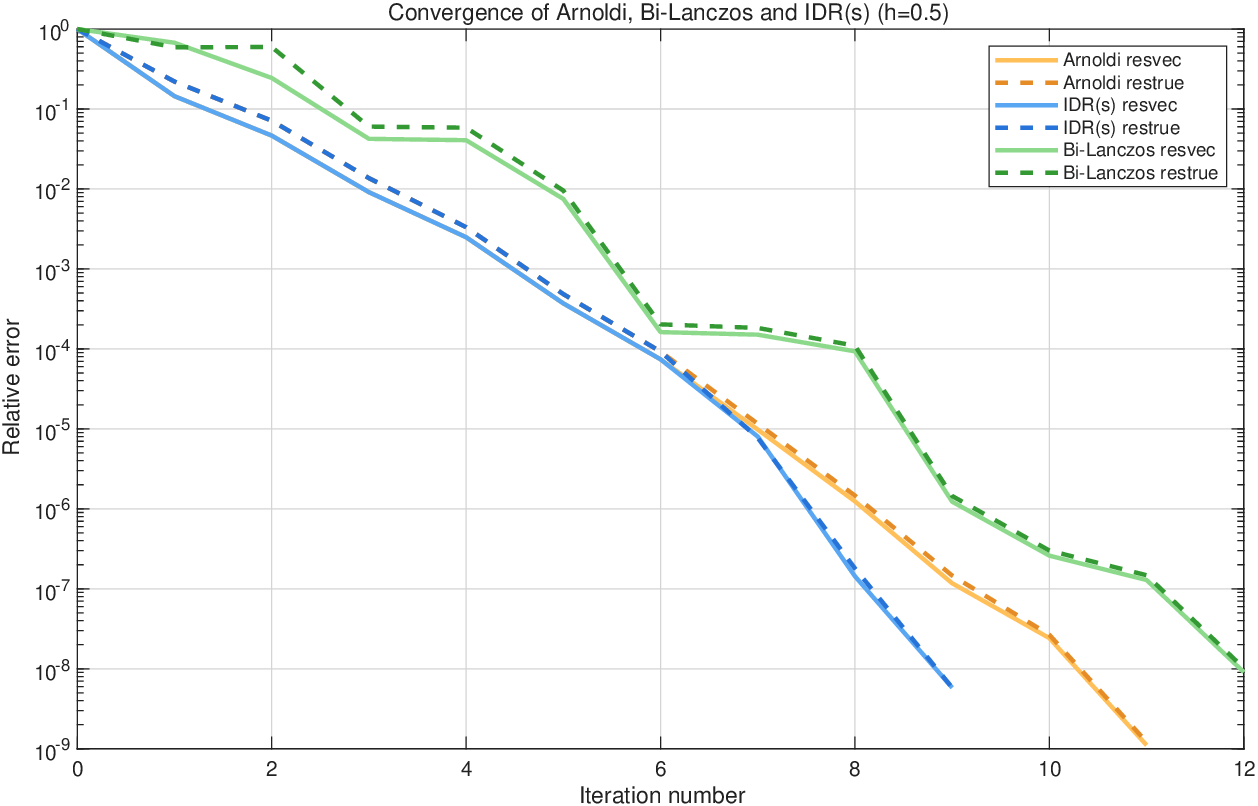}
        \caption*{\tt $h=0.5$}
    \end{subfigure}\hfill
    \begin{subfigure}[b]{0.31\textwidth}
        \centering
        \includegraphics[width=\textwidth]{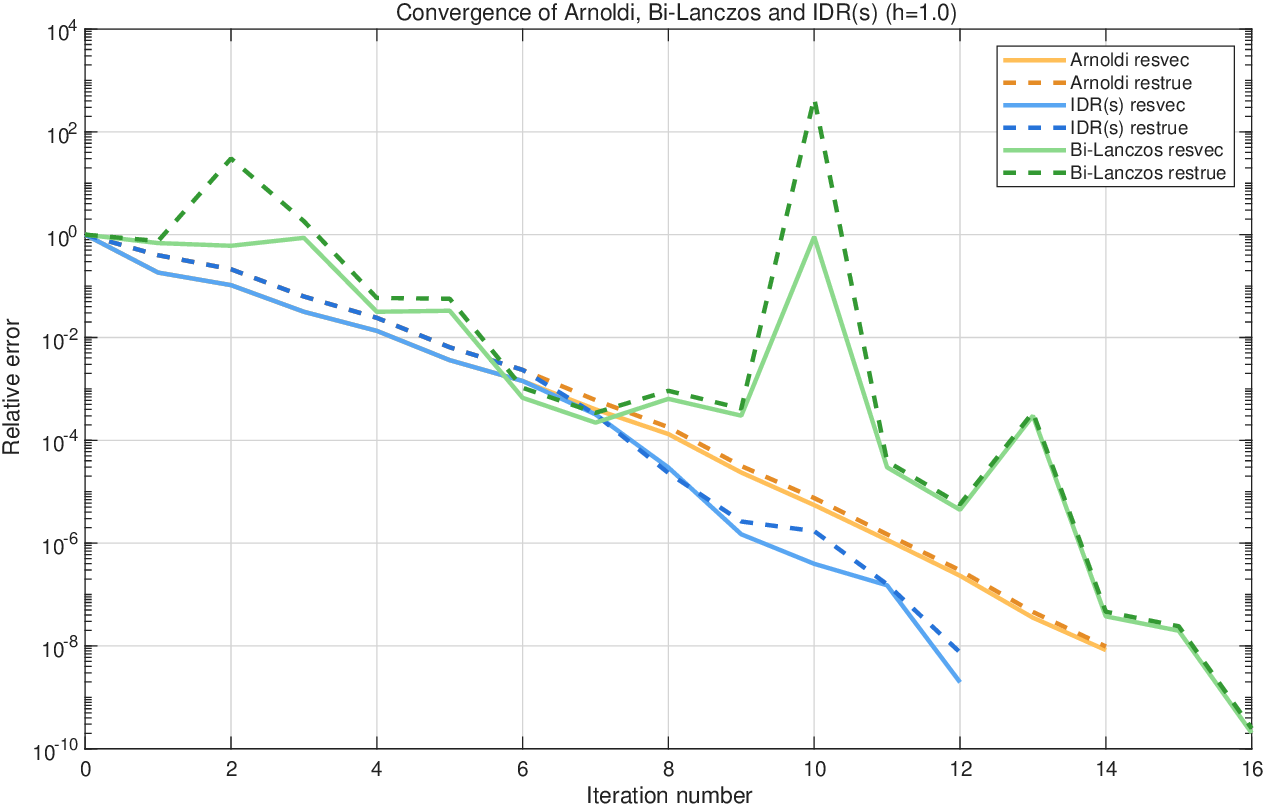}
        \caption*{\tt $h=1$}
    \end{subfigure}
    
    \vspace{1em}
    
    \begin{subfigure}[b]{0.31\textwidth}
        \centering
        \includegraphics[width=\textwidth]{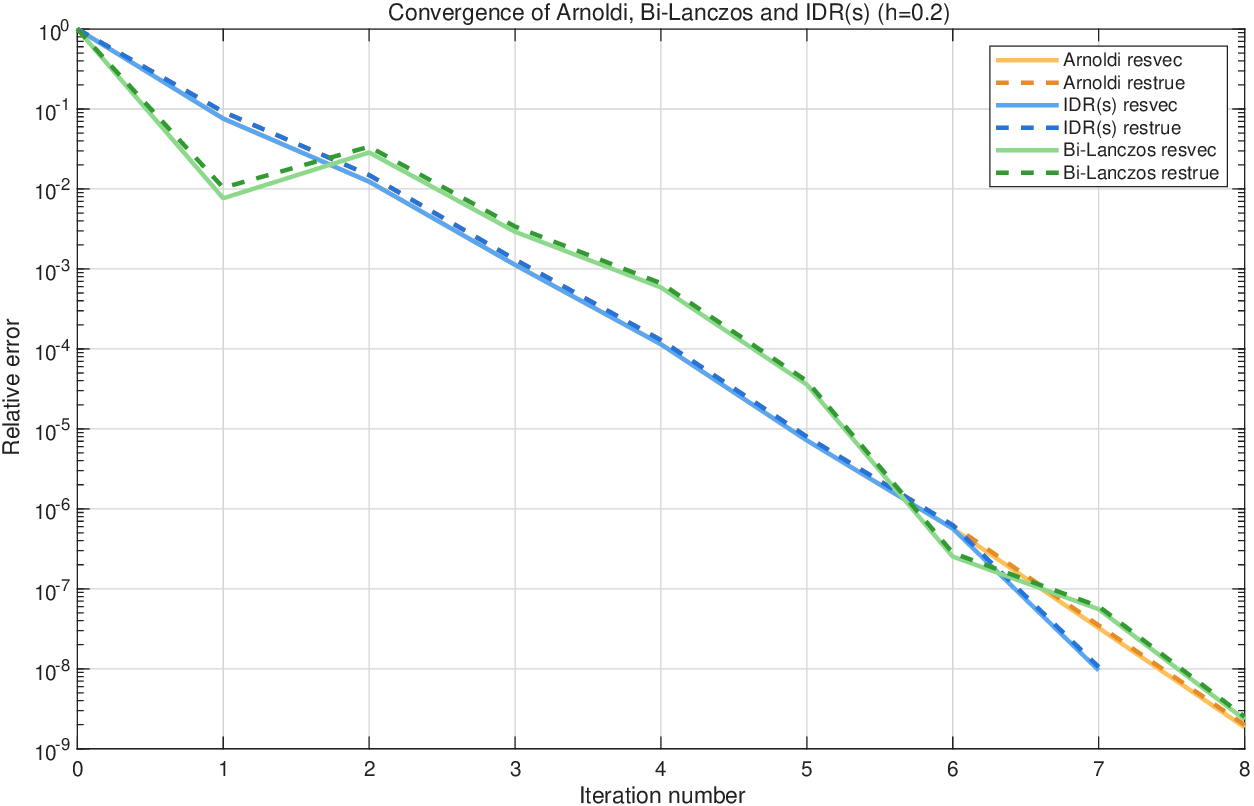}
        \caption*{\tt $h=0.2$}
    \end{subfigure}\hfill
    \begin{subfigure}[b]{0.31\textwidth}
        \centering
        \includegraphics[width=\textwidth]{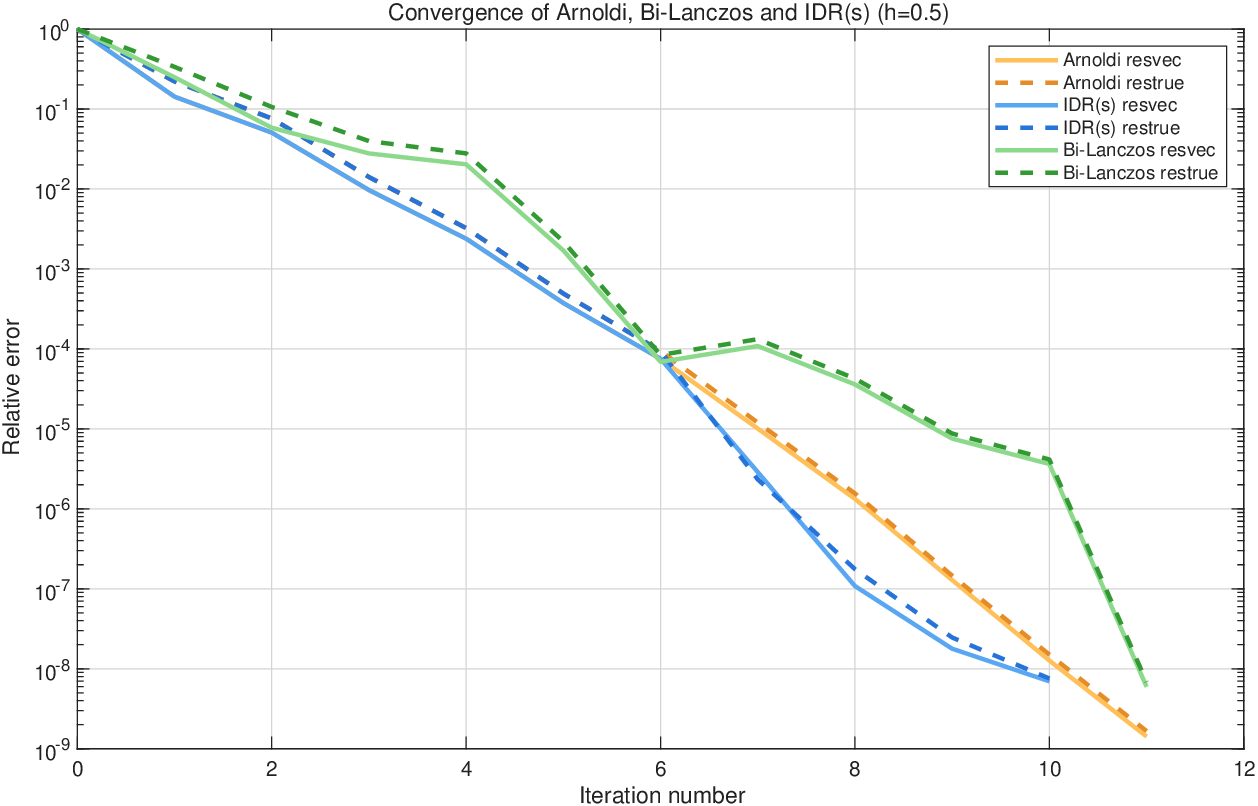}
        \caption*{\tt $h=0.5$}
    \end{subfigure}\hfill
    \begin{subfigure}[b]{0.31\textwidth}
        \centering
        \includegraphics[width=\textwidth]{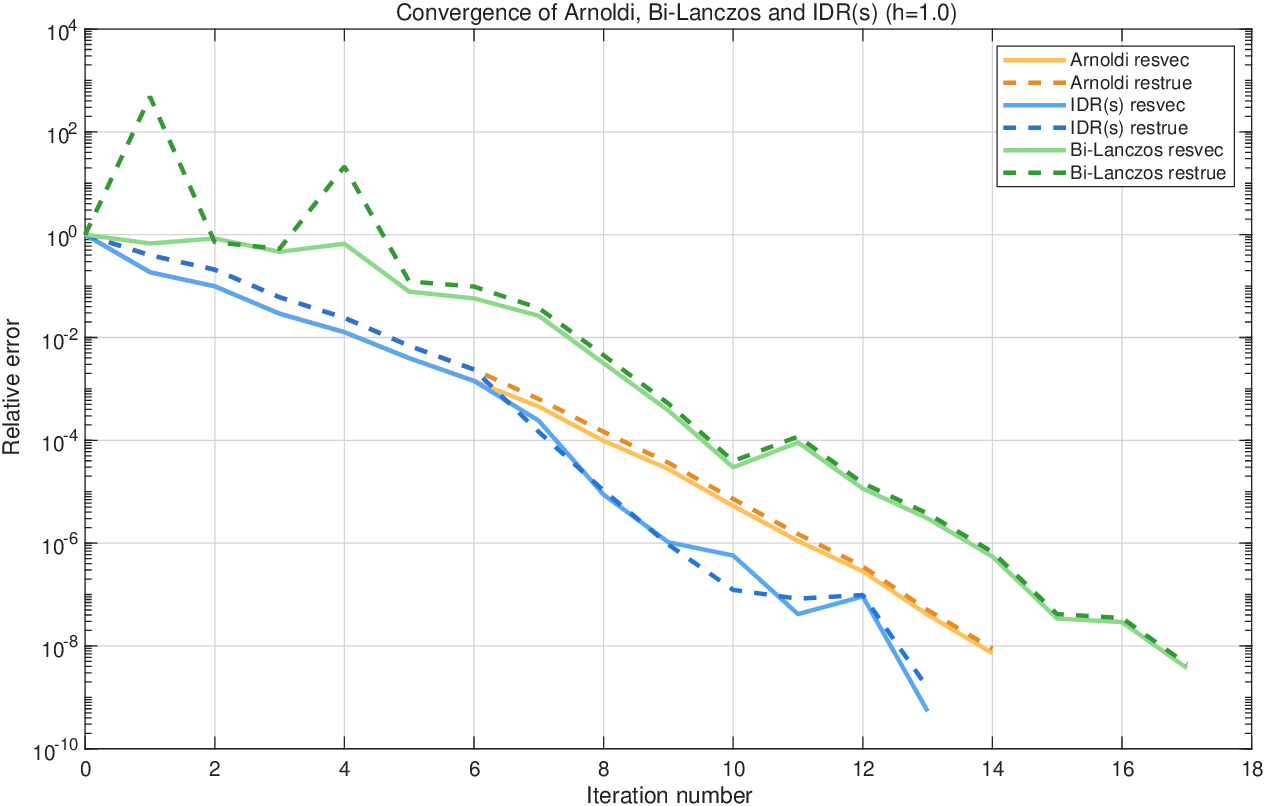}
        \caption*{\tt $h=1$}
    \end{subfigure}
    
    \vspace{1em}
    
    \begin{subfigure}[b]{0.31\textwidth}
        \centering
        \includegraphics[width=\textwidth]{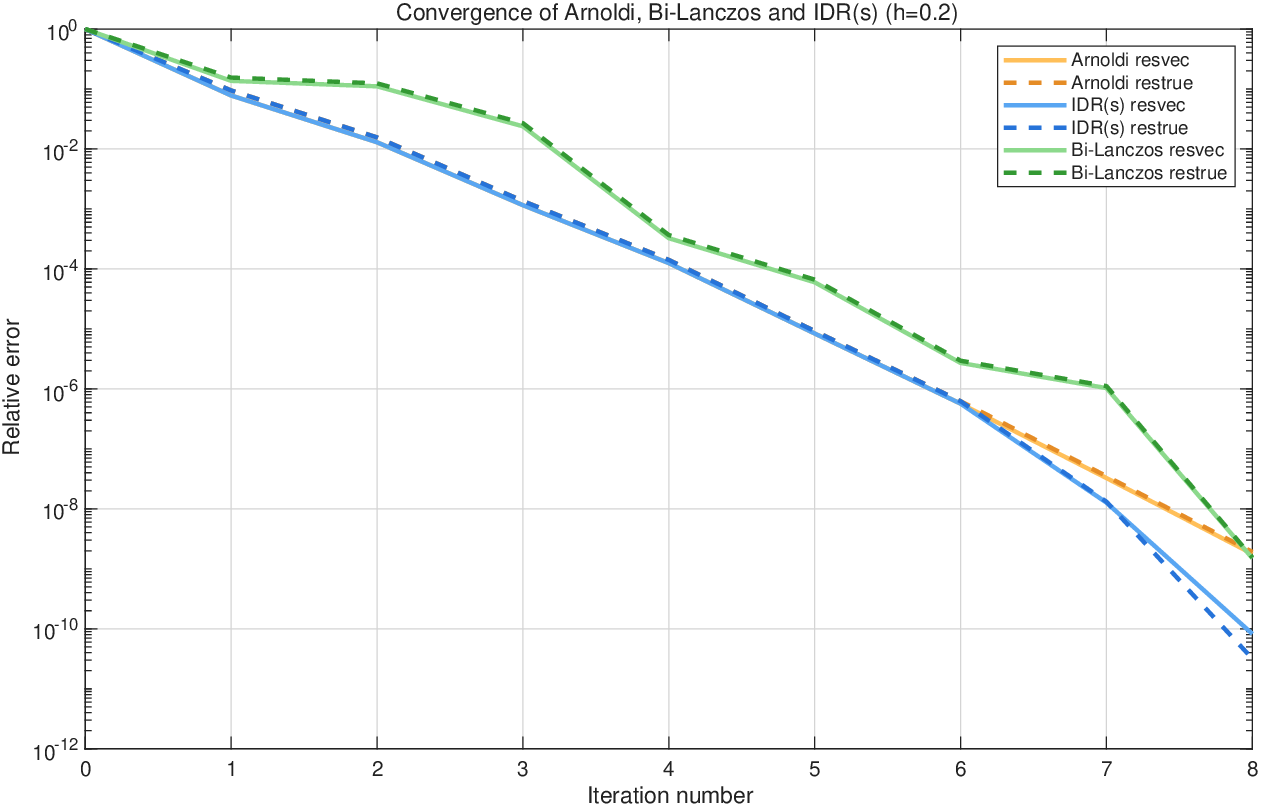} 
        \caption*{\tt $h=0.2$}
    \end{subfigure}\hfill
    \begin{subfigure}[b]{0.31\textwidth}
        \centering
        \includegraphics[width=\textwidth]{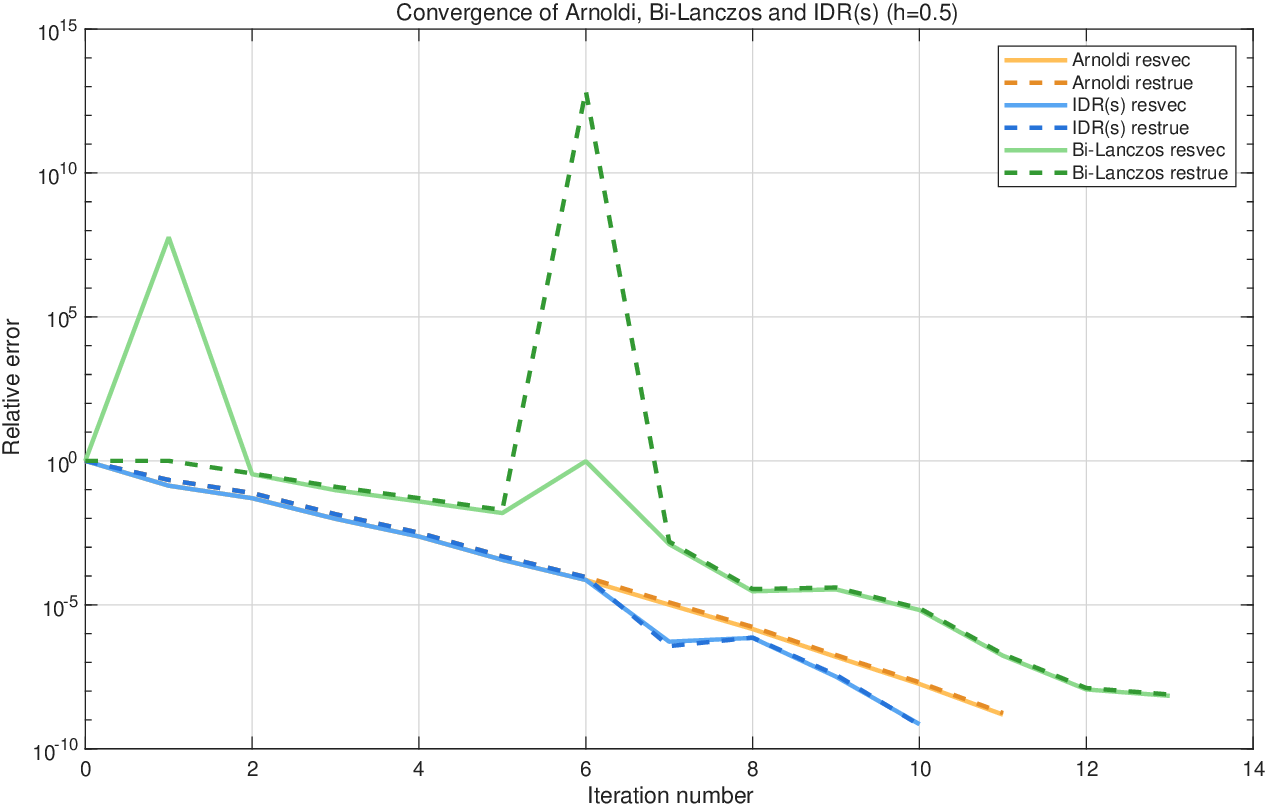}
        \caption*{\tt $h=0.5$}
    \end{subfigure}\hfill
    \begin{subfigure}[b]{0.31\textwidth}
        \centering
        \includegraphics[width=\textwidth]{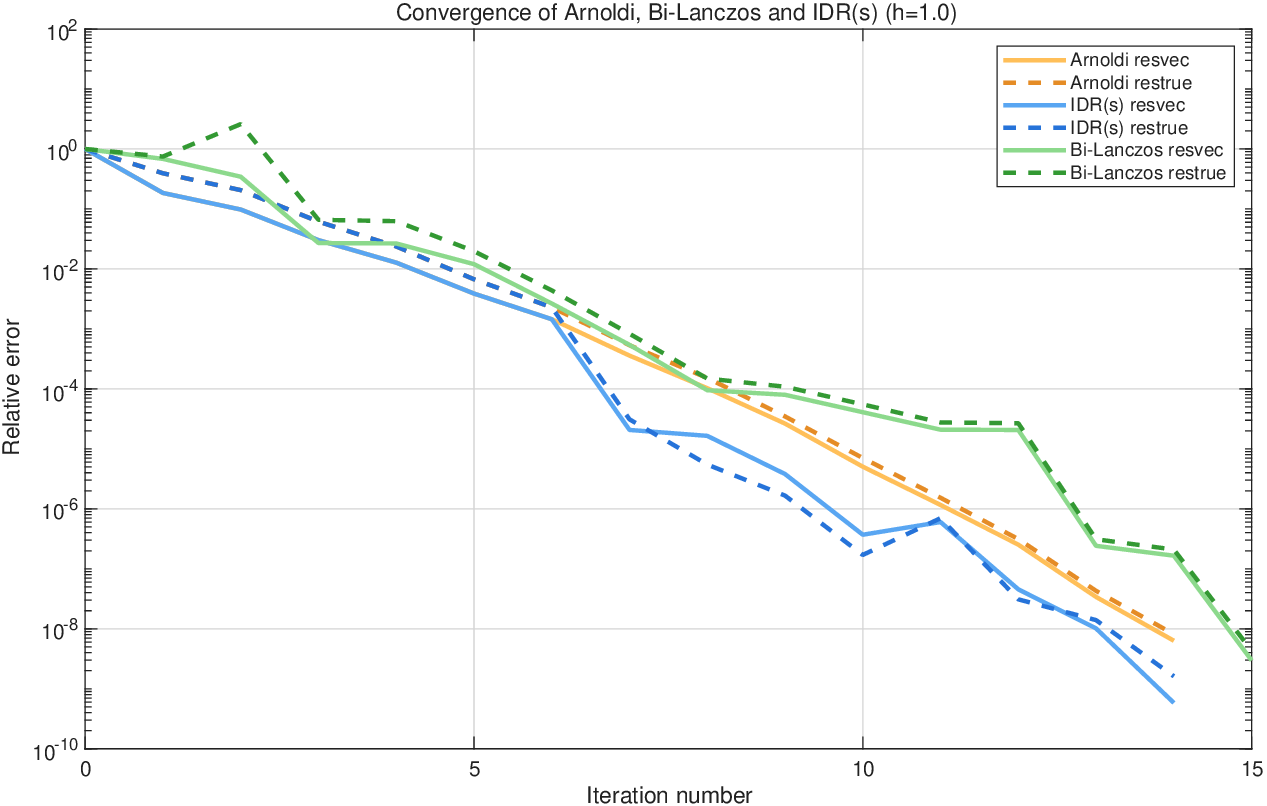}
        \caption*{\tt $h=1$}
    \end{subfigure}

    \caption{Convergence histories of three numerical methods for Example 1: \texttt{n = 10000} (top), \texttt{n = 15000} (middle), and \texttt{n = 20000} (bottom).}
    \label{fig1} 
\end{figure}

Fig. \ref{fig1} compares the convergence behaviors of the approximations based on the Arnoldi, IDR($s$) and bi-Lanczos processes, respectively. The bi-Lanczos method exhibits pronounced numerical instability, whereas the IDR($s$) a posteriori error estimate accurately tracks the true relative error.

\begin{table}[!htbp]
\centering
\caption{Three methods approximating $\mathbf{u}^\top \exp(hA)\mathbf{v}$ under varying $h$ for Example 1.}
\label{table1}
\begin{tabular}{ll cccccc}
\toprule 
\multirow{2}{*}{Size} & \multirow{2}{*}{Step size $h$} & \multicolumn{2}{c}{Arnoldi} & \multicolumn{2}{c}{IDR($s$)} & \multicolumn{2}{c}{Bi-Lanczos} \\
\cmidrule(lr){3-4} \cmidrule(lr){5-6} \cmidrule(lr){7-8} 
& & Iter & CPU (s) & Iter & CPU (s) & Iter & CPU (s) \\
\midrule 
\multirow{3}{*}{\texttt{n=10000}} 
& $h = 0.2$   & 8  & 0.2106 & 7  & 0.1619 & 9 & 0.3770 \\
& $h = 0.5$   & 11 & 0.2781 & 9 & 0.2308 & 12 & 0.4824 \\
& $h = 1$     & 14 & 0.3222 & 12 & 0.2572 & 16 & 0.6579\\
\midrule 
\multirow{3}{*}{\texttt{n=15000}} 
& $h = 0.2$   & 8  & 0.4167 & 7  & 0.3874 & 8 & 0.7109 \\
& $h = 0.5$   & 11 & 0.6576 & 10 & 0.5184 & 11 & 0.9851 \\
& $h = 1$     & 14 & 0.8594 & 13 & 0.6141 & 17 & 1.4949 \\
\midrule 
\multirow{3}{*}{\texttt{n=20000}} 
& $h = 0.2$   & 8  & 0.7656 & 8  & 0.6854 & 8 & 1.2941 \\
& $h = 0.5$   & 11 & 1.1928 & 10 & 0.8818 & 13 & 2.2500 \\
& $h = 1$     & 14 & 1.4236 & 14 & 1.0425 & 15 & 2.4919 \\
\bottomrule 
\end{tabular}
\end{table}

\textbf{Example 2}. This example evaluates the efficiency of Krylov subspace approximations for the matrix function bilinear form $\mathbf{u}^{\top}\exp(A)\mathbf{v}$, which dictates the option pricing dynamics within the Heston model. Here, the right vector $\mathbf{v} \in \mathbb{R}^n$ represents the discretized terminal payoff, while the sparse left vector $\mathbf{u} \in \mathbb{R}^n$ acts as a selector extracting the option value at the initial state. The scaled operator $A = \frac{1}{4}GT \in \mathbb{R}^{n \times n}$ is derived from a two-dimensional CTMC discretization \cite{Cui2021}, where $G$ denotes the transition generator matrix and $T$ is the time to maturity. This matrix $A$ exhibits highly specific structural characteristics: it is exceedingly sparse, {\em non-symmetric}, and block-tridiagonal. By fixing the variance grid size at $M=90$ and dynamically scaling the asset grid $N=70$, the global system dimension reaches $n = M \times N = 6,300$. 
\begin{figure}[!htp]
    \centering
    \includegraphics[width=0.40\textwidth]{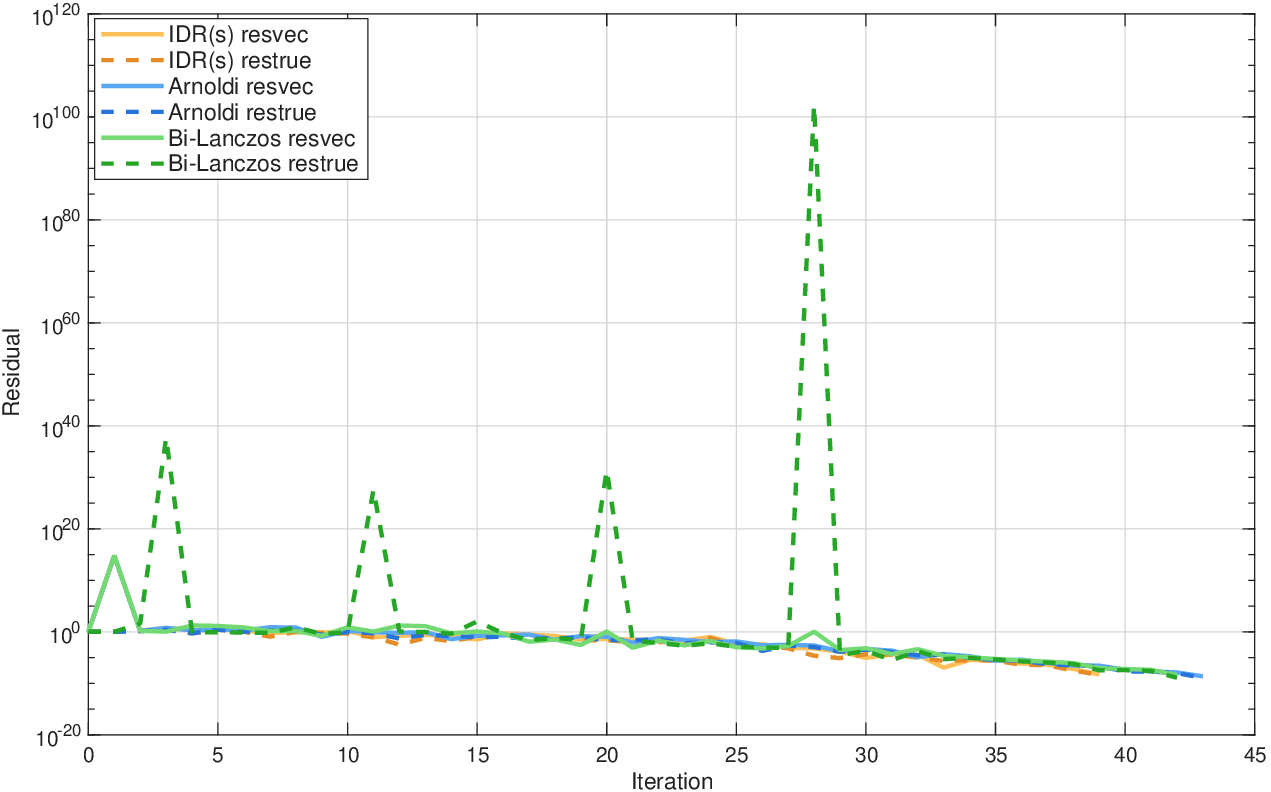} 
    \includegraphics[width=0.40\textwidth]{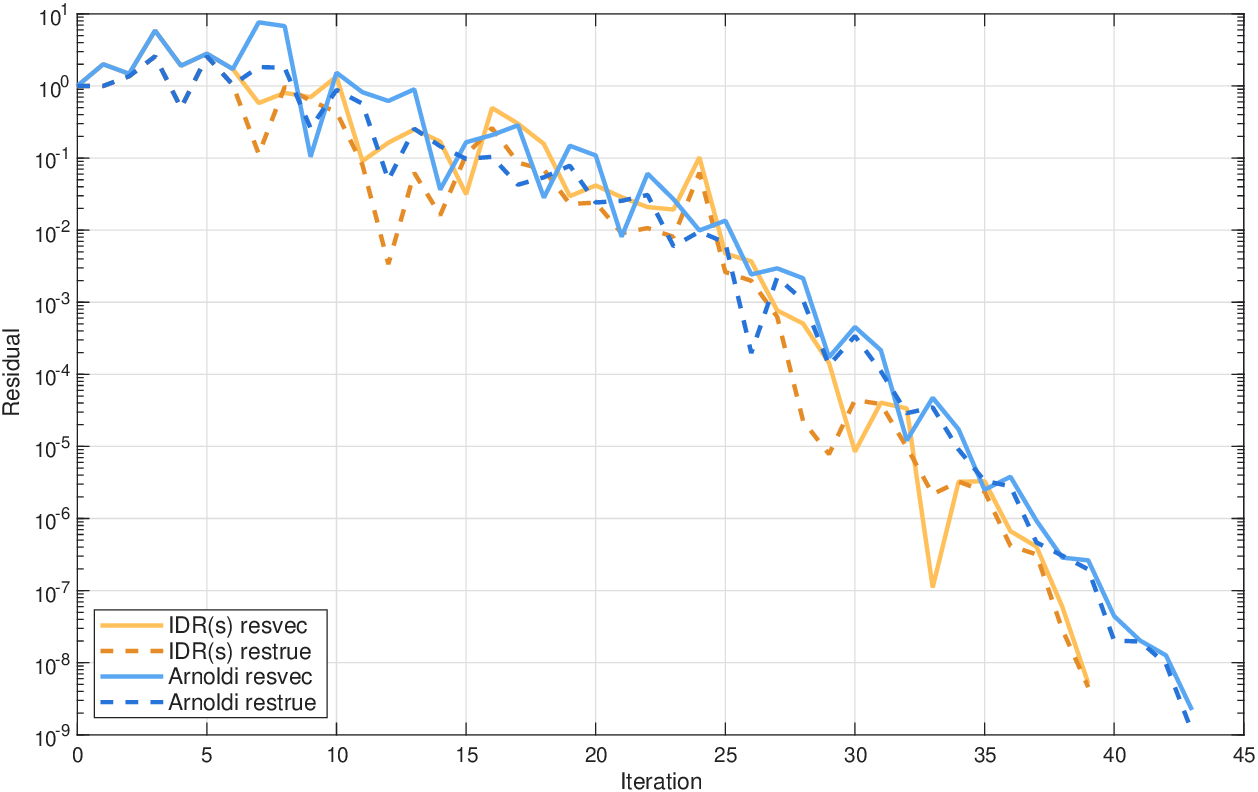} 
    \caption{Convergence histories of three numerical methods for Example 2.}
    \label{fig2}
\end{figure}
\begin{table}[!htbp]
\centering
\caption{Three methods approximating $\mathbf{u}^\top \exp(A)\mathbf{v}$ for Example 2.}
\label{table2}
\begin{tabular}{l cccccc}
\toprule 
\multirow{2}{*}{Computed Value} & \multicolumn{2}{c}{Arnoldi} & \multicolumn{2}{c}{IDR($s$)} & \multicolumn{2}{c}{Bi-Lanczos} \\
\cmidrule(lr){2-3} \cmidrule(lr){4-5} \cmidrule(lr){6-7} 
& Iter & CPU (s) & Iter & CPU (s) & Iter & CPU (s) \\
\midrule 
7.8577 & 43 & 0.1042 & 38 & 0.0840 & 42 & 0.0777 \\
\bottomrule 
\end{tabular}
\end{table}

Fig. \ref{fig2} and Table \ref{table2} compare the convergence and computational metrics of the approximations based on the Arnoldi, IDR($s$), and bi-Lanczos methods. All solvers successfully converge to the identical approximation. However, while Arnoldi and IDR($s$) converge smoothly, bi-Lanczos suffers from severe numerical instability. IDR($s$) proves to be the most efficient and reliable method; its recursive residuals accurately track true errors, and it converges in the fewest iterations. Consequently, the IDR($s$)-based method is the most efficient among the tested methods in this experiment.

\textbf{Example 3}. The test matrix employed in this section is identical to the one used in Example 1, with the dimension fixed at $n = 15000$. However, we now consider the approximation of the bilinear form for the matrix cosine, namely $\mathbf{u}^\top \cos(hG) \mathbf{v}$. To comprehensively evaluate the numerical stability and convergence behavior of the algorithms, the scaling parameter $h$ is systematically chosen as \(h = 0.2, 0.5, 1\). 
\begin{figure}[!htp]
    \centering
    \begin{minipage}[b]{0.3\textwidth}        \includegraphics[width=\textwidth]{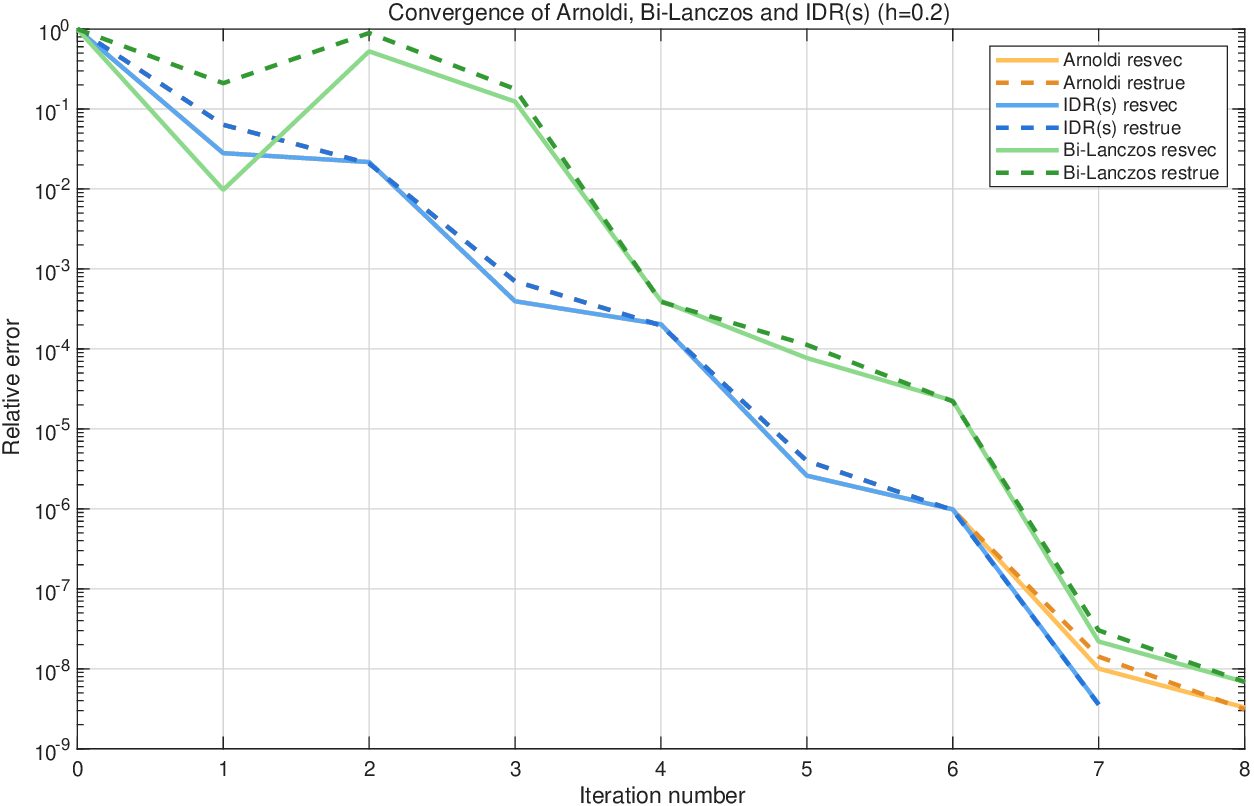}
        \centering{\tt $h=0.2$}
    \end{minipage}\hfill
    \begin{minipage}[b]{0.3\textwidth}
        \includegraphics[width=\textwidth]{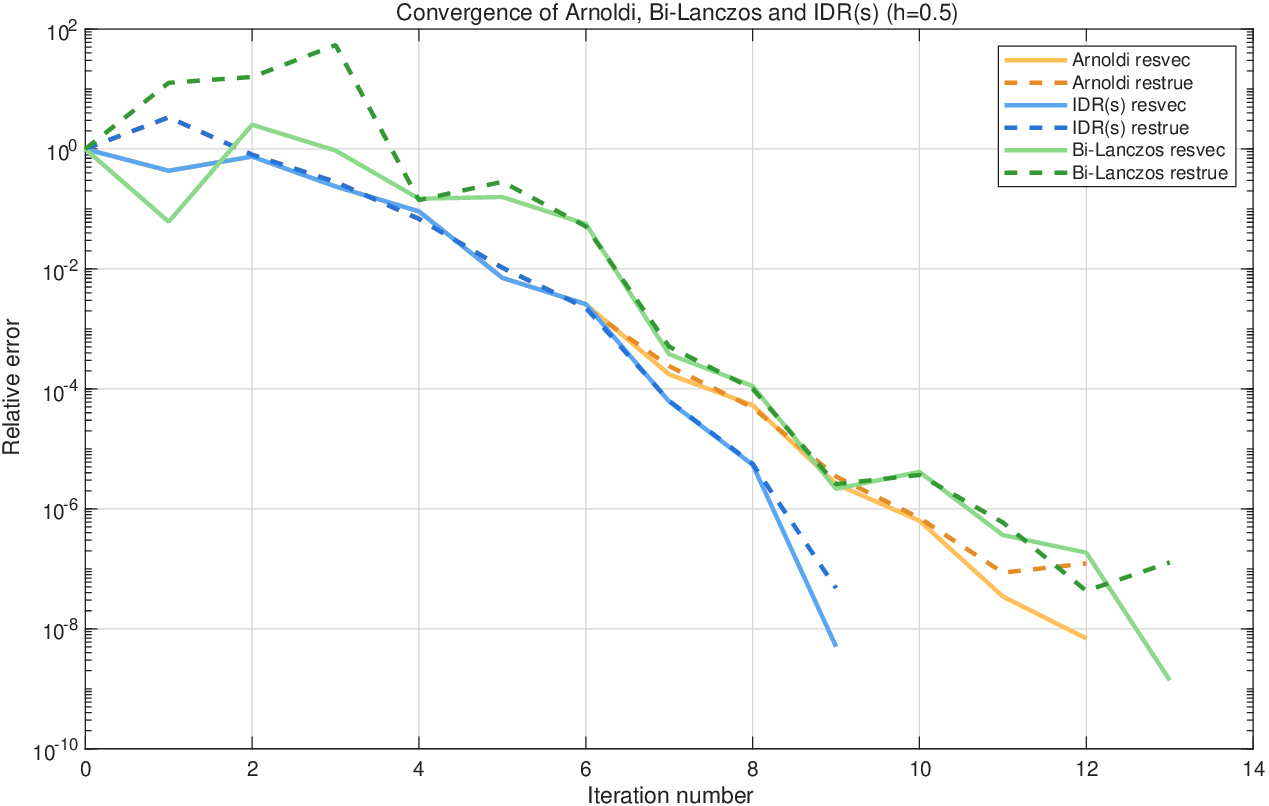}
        \centering{\tt $h=0.5$}
    \end{minipage}\hfill
    \begin{minipage}[b]{0.3\textwidth}
        \includegraphics[width=\textwidth]{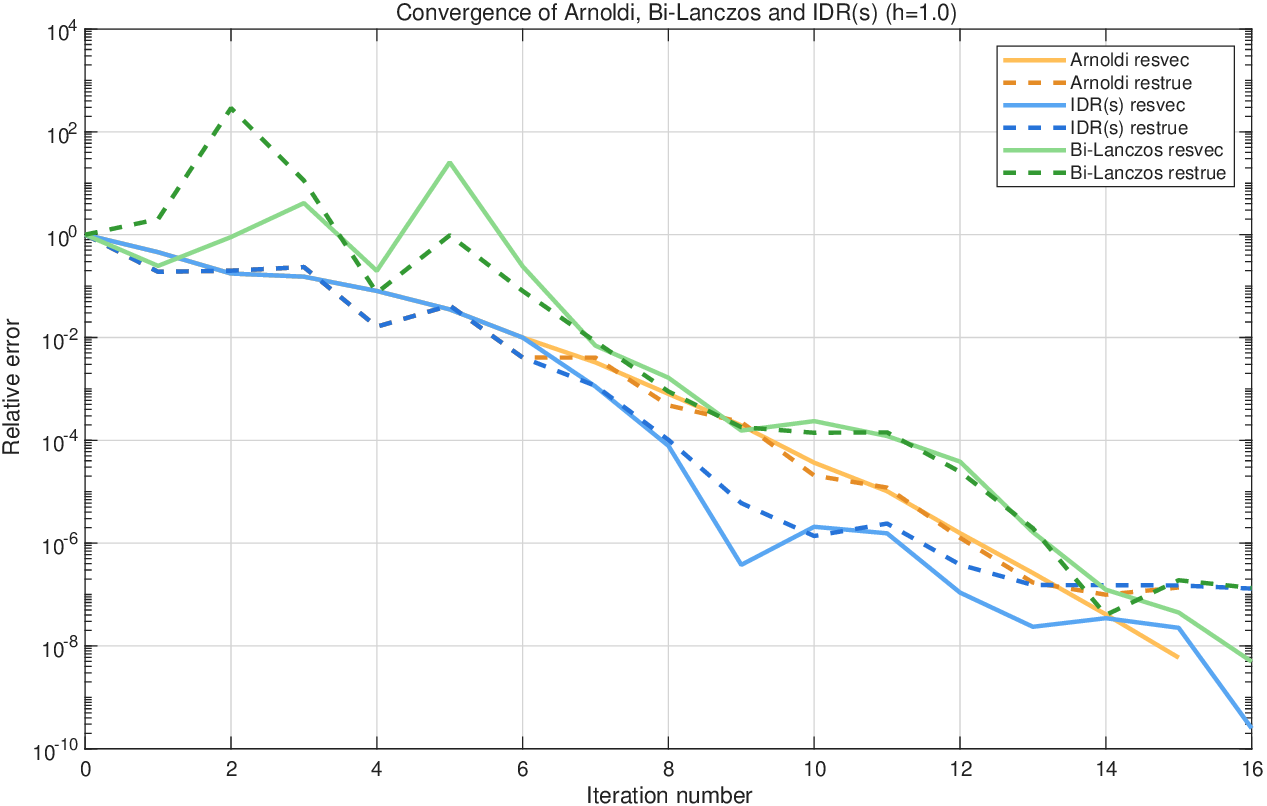}
        \centering{\tt $h=1$}
    \end{minipage}
    \caption{Convergence histories of three numerical methods for Example 3.}
    \label{fig3-1}
\end{figure}
\begin{table}[!htbp]
\centering
\setlength{\tabcolsep}{4.5mm}
\caption{Three methods approximating $\mathbf{u}^\top \cos(hA)\mathbf{v}$ under varying $h$ for Example 3.}
\label{table3}
\begin{tabular}{l cccccc}
\toprule 
\multirow{2}{*}{Step size $h$} & \multicolumn{2}{c}{Arnoldi} & \multicolumn{2}{c}{IDR($s$)} & \multicolumn{2}{c}{Bi-Lanczos} \\
\cmidrule(lr){2-3} \cmidrule(lr){4-5} \cmidrule(lr){6-7} 
& Iter & CPU (s) & Iter & CPU (s) & Iter & CPU (s) \\
\midrule 
$h=0.2$  & 8  & 0.6664 & 7  & 0.3844 & 8  & 0.8100 \\
$h=0.5$  & 12 & 1.9313 & 9  & 1.1522 & 13 & 2.3931 \\
$h=1$    & 15 & 1.6097 & 16 & 0.8547 & 16 & 1.4796 \\
\bottomrule 
\end{tabular}
\end{table}
Fig. \ref{fig3-1} and Table \ref{table3} compare the different approximations for $\mathbf{u}^\top \cos(hG)\mathbf{v}$. Unlike the oscillatory and inefficient bi-Lanczos method, both Arnoldi and IDR($s$) converge stably. Furthermore, IDR($s$) features a reliable a posteriori error estimate and consistently requires the least CPU time across all step sizes. Ultimately, the IDR($s$)-based method shows the best overall efficiency among the tested methods.

\textbf{Example 4}. In this experiment, we evaluate our method on a sequence of nonsymmetric test matrices (\texttt{cage13}, \texttt{cage14} and \texttt{cage15}) from the \texttt{vanHeukelum} group in the SuiteSparse Matrix Collection \footnote{\url{https://sparse.tamu.edu/vanHeukelum}}. The projection vectors $\mathbf{u}, \mathbf{v} \in \mathbb{R}^n$ are randomly generated and normalized such that $\|\mathbf{u}\|_2 = \|\mathbf{v}\|_2 = 1$. Due to the massive dimensions of the \texttt{cage} matrix series, computing $\mathbf{u}^\top \sqrt{A} \mathbf{v}$ is computationally intractable. To establish a high-precision reference solution for our convergence analysis, we employ a numerical contour integration approach based on Gaussian quadrature. Specifically, the reference value is computed using a Gauss-Legendre quadrature rule, where the nodes and weights are generated via the Golub-Welsch algorithm and algebraically mapped from the standard interval $[-1, 1]$ to the semi-infinite interval $[0, \infty)$. This procedure elegantly decouples the matrix function evaluation into a series of independent shifted linear systems, which are subsequently solved to high precision \cite{higham2008}. Furthermore, for the evaluation of the analytic function sequence $\phi_i(H_m)$ required in the a posteriori error estimate, the expansion point $t_0$ is adaptively chosen as $h_{1,1}$.

\begin{figure}[!htp]
    \centering
    \begin{minipage}[b]{0.3\textwidth}
        \includegraphics[width=\textwidth]{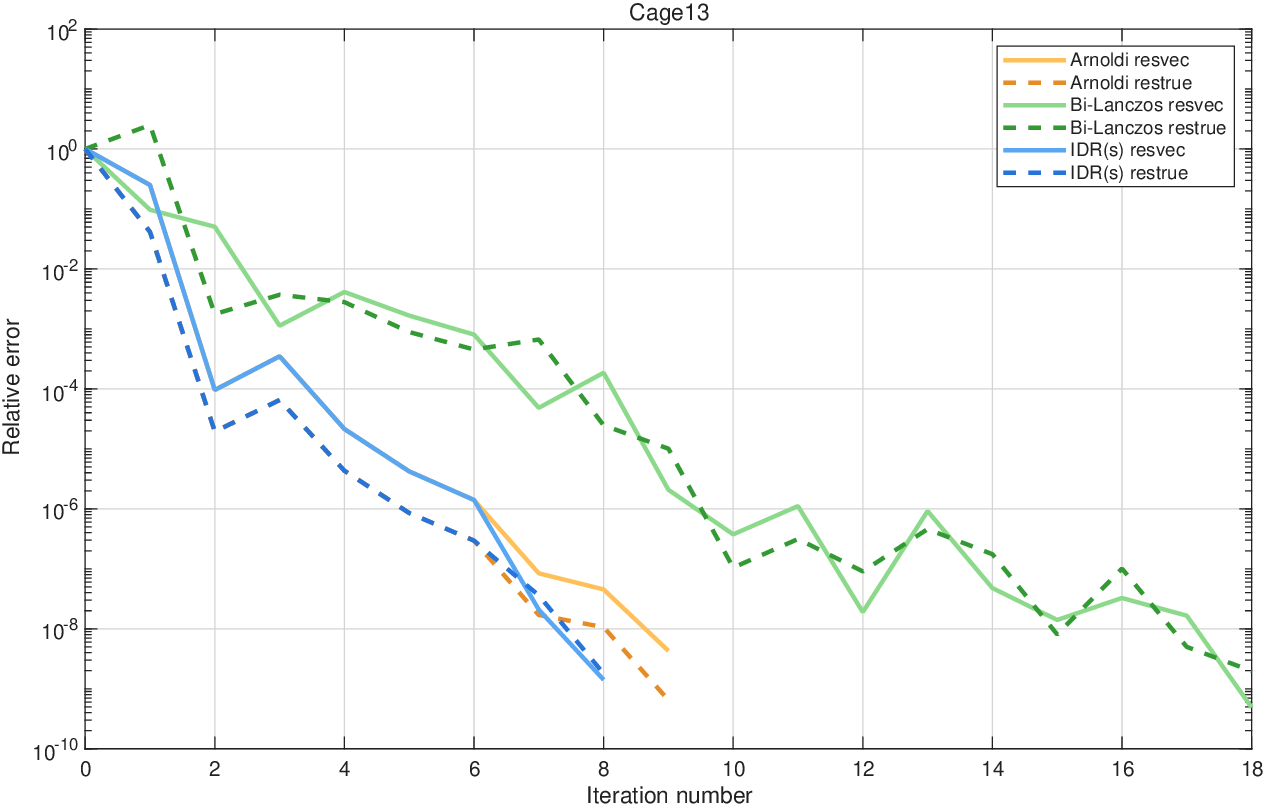}
        \centering{\tt $cage13$}
   \end{minipage}\hfill
    \begin{minipage}[b]{0.3\textwidth}
        \includegraphics[width=\textwidth]{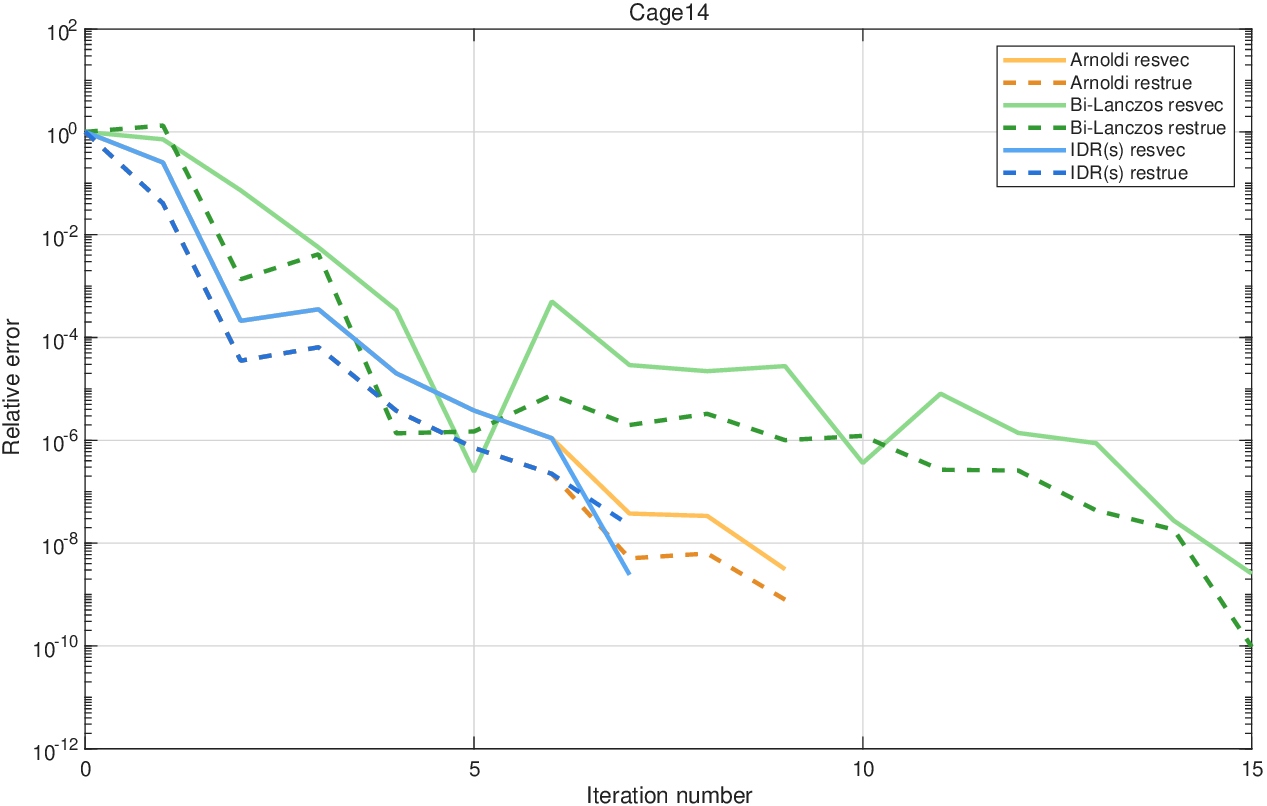}
        \centering{\tt $cage14$}
    \end{minipage}\hfill
    \begin{minipage}[b]{0.3\textwidth}
        \includegraphics[width=\textwidth]{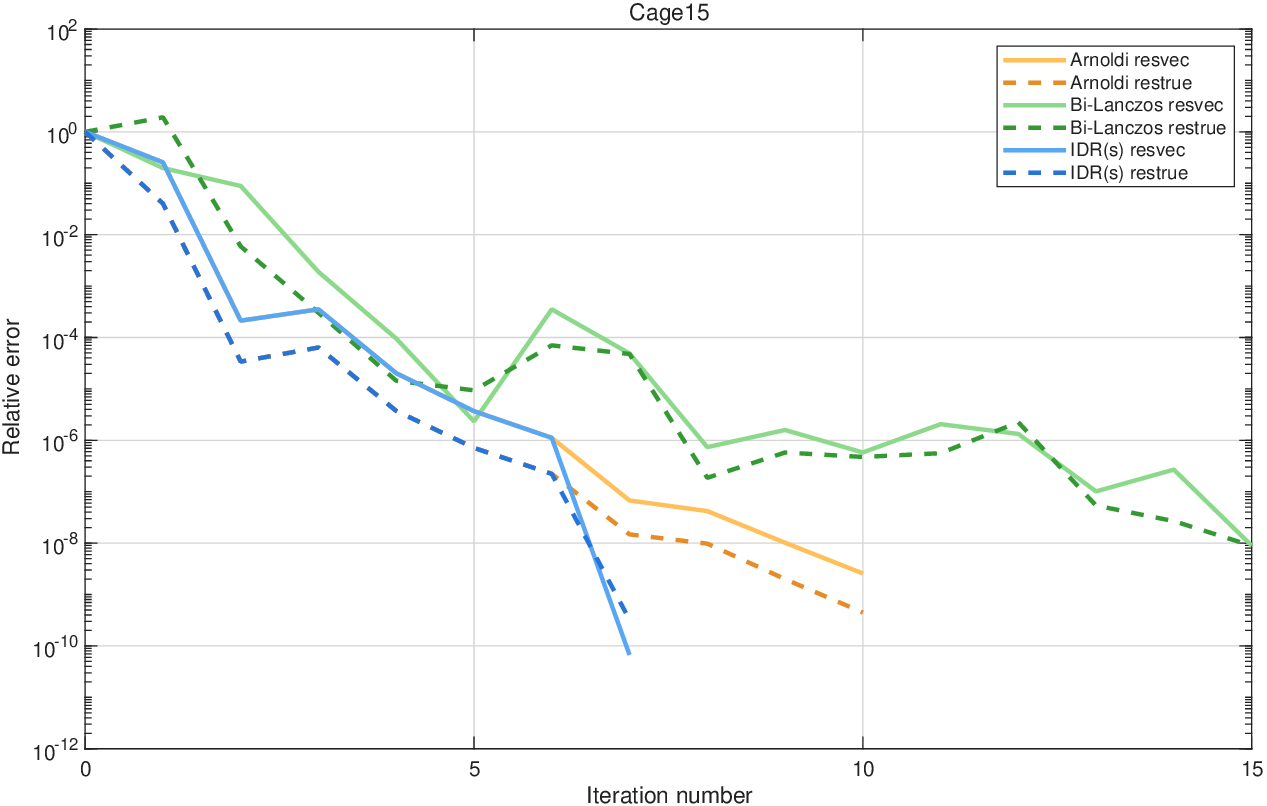}
        \centering{\tt $cage15$}
    \end{minipage}
    \caption{Convergence histories of three numerical methods for Example 4.}
    \label{fig4-1}
\end{figure}
\begin{table}[!htbp]
\centering
\setlength{\tabcolsep}{4.5mm}      
\caption{Three methods approximating $\mathbf{u}^\top \sqrt{A}\mathbf{v}$ for Example 4.} 
\label{table4}
\begin{tabular}{l cccccc}
\toprule 
\multirow{2}{*}{Matrix} & \multicolumn{2}{c}{Arnoldi} & \multicolumn{2}{c}{IDR($s$)} & \multicolumn{2}{c}{Bi-Lanczos} \\
\cmidrule(lr){2-3} \cmidrule(lr){4-5} \cmidrule(lr){6-7} 
& Iter & CPU (s) & Iter & CPU (s) & Iter & CPU (s) \\
\midrule 
\texttt{cage13} & 9  & 0.2174 & 8 & 0.1893 & 18 & 0.9602 \\
\texttt{cage14} & 9  & 0.6266 & 7 & 0.4566 & 15 & 1.9539 \\
\texttt{cage15} & 10 & 2.6250 & 7 & 1.6627 & 15 & 6.9037 \\
\bottomrule 
\end{tabular}
\end{table}
Numerical results show that the Gauss quadrature method takes much more CPU time than the approximations based on the  Arnoldi, IDR($s$), and bi-Lanczos methods. For the Krylov subspace methods (see Fig. \ref{fig4-1} and Table \ref{table4}), Bi-Lanczos is slow and oscillatory. Conversely, the IDR($s$)-based method converges smoothly with reliable error bounds, consistently requiring less CPU time than Arnoldi, making it particularly effective for the tested large-scale matrices.

Despite the absence of theoretical guarantees on the boundedness of higher-order derivatives or the formal reliability of the leading-order error term, our numerical experiments confirm that the resulting error expansion serves as a viable and effective convergence heuristic for $\mathbf{u}^\top \sqrt{A}\mathbf{v}$.

\textbf{Example 5}. We construct a non-symmetric block-diagonal matrix \(A=\operatorname{blkdiag}(A_{\rm act},B)\) of order 5000. The active block \(A_{\rm act}=QJQ^\top\in\mathbb{R}^{80\times80}\) is a randomly rotated scaled nilpotent Jordan block with subdiagonal entries \(\alpha=2\), which preserves a controllable Krylov-chain structure. The background block \(B\) is a variable-coefficient banded upper-Hessenberg matrix, making the full matrix large-scale and nonnormal. The left and right vectors \(\mathbf{u}\) and \(\mathbf{v}\) are designed to induce a controlled temporary cancellation in the error estimator at the sixth iteration, while a nonzero tail in \(\mathbf{u}\) keeps the true relative error non-negligible. 

\begin{table}[!htbp]
\centering
\setlength{\tabcolsep}{4.0mm} 
\caption{Comparison of single and dual stopping criteria for Example 5.}
\label{table5}
\begin{tabular}{ll cccc} 
\toprule 
\multirow{2}{*}{Method} & \multirow{2}{*}{Criterion} & \multirow{2}{*}{Iter} & \multirow{2}{*}{CPU (s)} & \multicolumn{2}{c}{Relative error} \\
\cmidrule(lr){5-6} 
& & & & Estimated & True \\
\midrule 
\multirow{2}{*}{Arnoldi} 
& Single & 6  & 0.0047 & 1.4799e-17 & 4.1257e-02 \\
& Dual   & 17 & 0.0279 & 9.7104e-11 & 1.0975e-10 \\
\midrule 
\multirow{2}{*}{IDR($s$)} 
& Single & 6  & 0.0052 & 1.4799e-17 & 4.1257e-02 \\
& Dual   & 17 & 0.0144 & 1.3733e-09 & 1.9373e-09 \\
\midrule 
\multirow{2}{*}{Bi-Lanczos} 
& Single & 17 & 0.0192 & 4.5777e-09 & 5.3221e-09 \\
& Dual   & 18 & 0.0194 & 1.0365e-10 & 1.2029e-10 \\
\bottomrule 
\end{tabular}
\end{table}

\begin{figure}[!htp]
    \centering
    
    \includegraphics[width=0.80\textwidth]{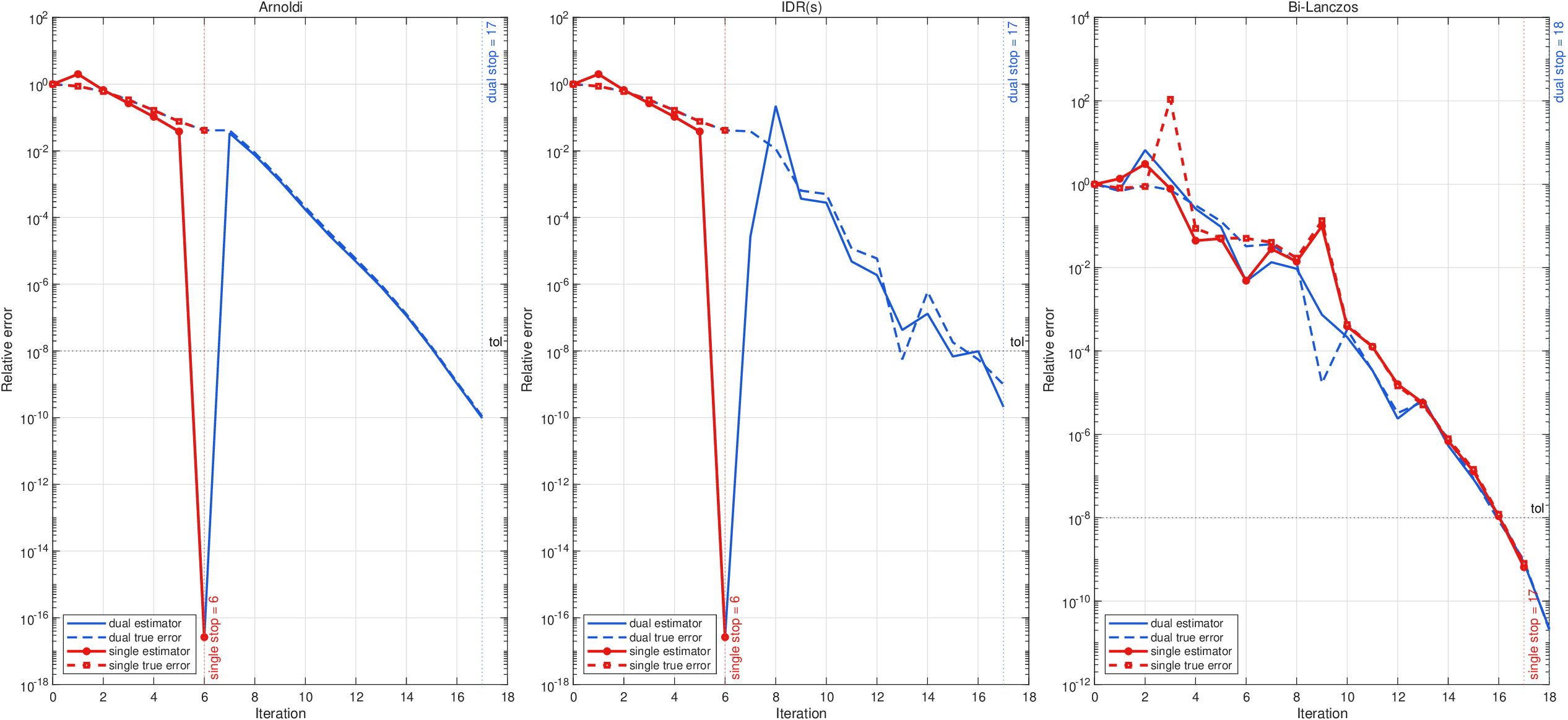} 
    \caption{Convergence histories of three numerical methods under single and dual stopping criteria.}
    \label{fig5}
\end{figure}

The stopping behavior of the three methods is summarized in Table \ref{table5} and Fig. \ref{fig5}. For the Arnoldi and IDR($s$) methods, the single stopping criterion produces false convergence at iteration 6 because of a spuriously small estimated error of approximately \(10^{-17}\), whereas the true relative error remains as large as \(4.13\times10^{-2}\). In contrast, the dual criterion avoids premature termination caused by this temporary cancellation and continues the iterations until the true relative error falls below the prescribed tolerance of \(10^{-8}\). For bi-Lanczos, both criteria achieve the prescribed accuracy, although the dual criterion further reduces the final true relative error. 
\section{Conclusions}
\label{sec:conclusions}
This paper studied the computation of bilinear forms of matrix functions using the IDR($s$) method and developed a corresponding a posteriori error indicator and stopping criterion. Based on the Hessenberg decomposition generated by the IDR($s$) process, we constructed a Krylov subspace approximation for \eqref{1.1} and derived an explicit error expansion for the resulting approximation. We proved that, under mild analyticity and boundedness assumptions, the remainder term of this expansion vanishes as the truncation order increases, thereby validating the expansion. Motivated by this analysis, the leading computable term was used as a practical a posteriori error indicator, and a dual stopping criterion was formulated to reduce the risk of premature termination caused by accidental cancellations. The numerical experiments confirmed that the proposed error indicator tracks the true error well for several representative matrix functions, including the exponential, trigonometric, and square-root functions. In particular, the controlled cancellation experiment showed that the dual stopping criterion avoids false convergence in cases where the single criterion terminates too early. These results demonstrate the practical usefulness of the proposed IDR($s$)-based framework for large-scale bilinear forms of matrix functions. As an oblique projection method, however, the IDR($s$) process does not preserve global orthogonality of the basis vectors, and highly ill-conditioned cases may still lead to numerical instability or parameter sensitivity. This limitation motivates further investigation of stabilization strategies for IDR($s$)-based matrix-function approximations. 

\section*{Acknowledgments}
We acknowledge Dr. Wensheng Yang for providing the MATLAB code for Example 2, which enabled us to generate the corresponding test matrices and the two vectors. 
\section*{Author Contributions}
Q.Xue, X.Yue and XM.Gu wrote and revised the main manuscript text; Q.Xue and XM.Gu prepared numerical experiments. All the authors reviewed the manuscript.
\section*{Funding}
The work is supported by the Science Challenge Project (No. TZ2024009), the National Natural
Science Foundation of China (Nos. 72431008 and 12371373), and the National Key R\&D Program of China (No. 2023YFB3001702).

\section*{Ethics declarations}
\subsection*{Conflict of interest}
The authors declare no competing interests.
\subsection*{Data availability}
The data that support the findings of this study are available in the article. The source code used in this study is available from the authors upon reasonable request.

\appendix
\section{Convergence Analysis of the Remainder Term} \label{app:convergence_proof}
In this appendix, we present the technical details of the convergence rate analysis for the remainder term $|\mathbf{u}^\top p_i(A)\mathbf{s}_m^i|$ in Theorem \ref{Theorem 1}, covering both entire and general analytic functions.

\subsection{Proof for functions satisfying the derivative-growth condition (e.g., $e^z$, $\cos z$)}\label{A.1}
We first establish the geometric convergence for entire functions. Since the matrix $A$ is non-symmetric and generally possesses complex eigenvalues, we utilize the complex Schur decompositions:
\begin{equation} \label{eq:complex_schur}
Q_1^* A Q_1 = U_1 \quad \text{and} \quad Q_2^* H_m Q_2 = U_2,
\end{equation}
where $Q_1 \in \mathbb{C}^{n \times n}$ and $Q_2 \in \mathbb{C}^{m \times m}$ are unitary matrices (i.e., $Q_1^* Q_1 = I_n$ and $Q_2^* Q_2 = I_m$), and $U_1, U_2$ are upper triangular matrices with eigenvalues of $A$ and $H_m$ on their diagonals, respectively. We denote $\hat{U}_1$ and $\hat{U}_2$ as their strictly upper triangular parts. 

By applying the derivative bounds \eqref{eq:phi_deriv_new} to the Schur-based matrix function bounds, we have:
\begin{align*}
\|\phi_i(A)\| &\le \sum_{k=0}^{n-1} M k! \rho^{-(k+i)} \frac{\|\hat{U}_1\|_F^k}{k!} = M \rho^{-i} \sum_{k=0}^{n-1} \left( \frac{\|\hat{U}_1\|_F}{\rho} \right)^k = \tilde{C}_1 \rho^{-i}, \\
\|\phi_i(H_m)\| &\le \sum_{k=0}^{m-1} M k! \rho^{-(k+i)} \frac{\|\hat{U}_2\|_F^k}{k!} = M \rho^{-i} \sum_{k=0}^{m-1} \left( \frac{\|\hat{U}_2\|_F}{\rho} \right)^k = \tilde{C}_2 \rho^{-i},
\end{align*}
where $\tilde{C}_1$ and $\tilde{C}_2$ are positive constants independent of $i$ (since the matrix sizes $n$ and $m$ are finite). 

From the expression of the error vector $\mathbf{s}_m^i = \phi_i(A)\mathbf{v} - \beta V_m \phi_i(H_m)\mathbf{e}_1$, and using the fact that $\|V_m\|_2 \le \sqrt{m}$ \cite{Astudillo2016}, it follows that:
\begin{align}
\|\mathbf{s}_m^i\| &\le \|\phi_i(A)\| \cdot \|\mathbf{v}\| + \beta \|V_m\|_2 \cdot \|\phi_i(H_m)\| \cdot \|\mathbf{e}_1\| \notag \\
&\le \tilde{C}_1 \rho^{-i} + \beta \sqrt{m} \tilde{C}_2 \rho^{-i} = \tilde{C}_3 \rho^{-i}. \label{eq:s_norm_bound_new}
\end{align}

Since the set of interpolation nodes $\{tj\}_{j=0}^{i-1}$ belongs to the closed convex set $S$, we can bound the polynomial term $\|p_i(A)\| = \|(A - t_0 I)\dots(A - t_{i-1}I)\| \le C_2^i$. Combining this with Eq. \eqref{eq:s_norm_bound_new} via the Cauchy-Schwarz inequality, the remainder term of the error series is bounded by:
\begin{equation}
\left| \mathbf{u}^\top p_i(A) \mathbf{s}_m^i \right| \le \|\mathbf{u}\|_2 \cdot \|p_i(A)\|_2 \cdot \|\mathbf{s}_m^i\| \le C_3 C_2^i \tilde{C}_3 \rho^{-i} = C_4 \left( \frac{C_2}{\rho} \right)^i.
\label{eqA3}
\end{equation}

Since we explicitly assume $\rho > C_2$ in Theorem \ref{Theorem 1}, the ratio satisfies $C_2 / \rho < 1$. Consequently, the remainder term vanishes geometrically to $0$ as $i \to \infty$:
\begin{equation*}
\lim_{i \to \infty} \left| \mathbf{u}^\top p_i(A) \mathbf{s}_m^i \right| = 0,
\end{equation*}
which rigorously completes the proof of geometric convergence for entire functions.

\subsection{Proof for General Analytic Functions (e.g., $z^{-1}$, $z^{1/2}$)}
\label{app:general_analytic_proof}
For general analytic functions with singularities or branch cuts, we assume $f(z)$ represents a chosen single-valued analytic branch. Let $r > 0$ denote the distance from $S$ to the boundary of the domain of analyticity of $f(z)$ (which contains any singularities or branch cuts). 

By the definition of the open neighborhood of analyticity, we define the $\rho$-neighborhood (or inflated neighborhood) of $S$ for any intermediate radius $0 < \rho < r$ as:
\begin{equation} \label{eq:S_rho_def}
S_\rho := \{ z \in \mathbb{C} : \text{dist}(z, S) \le \rho \}.
\end{equation}
Since $\rho < r$, the compact set $S_\rho$ is entirely contained within the domain of analyticity of $f(z)$. Therefore, Cauchy's integral formula can be strictly applied over the boundary $\partial S_\rho$, which forms a rectifiable closed contour enclosing $S$. This mathematically guarantees that the divided differences and their derivatives on $S$ exhibit geometric decay:
\begin{equation} \label{eq:cauchy_est_rho}
\sup_{t \in S} \left| \phi_i(t) \right| \le M_{\rho} \rho^{-i} \quad \text{and} \quad \sup_{z\in S} \left| \phi_i^{(k)}(z) \right| \le M_{k, \rho} \rho^{-i},
\end{equation}
where $M_{\rho}$ and $M_{k, \rho}$ are positive constants depending on $\rho$. 

Following a similar Schur decomposition bounding procedure using the complex Schur decompositions \eqref{eq:complex_schur} as detailed in Appendix \ref{A.1}, we obtain the bound for the error vector:
\begin{equation} \label{eq:s_rho_bound}
\|\mathbf{s}_m^i\| \le \tilde{C}_{\rho} \rho^{-i},
\end{equation}
where $\tilde{C}_{\rho}$ is a positive constant. 

By combining Eq. \eqref{eq:s_rho_bound} with the Cauchy-Schwarz inequality and the polynomial bound $\|p_i(A)\|_2 \le C_2^i$, the remainder term of the error series can be strictly bounded by:
\begin{equation} \label{A.4}
\left| \mathbf{u}^\top p_i(A) \mathbf{s}_m^i \right| \le \|\mathbf{u}\|_2 \cdot \|p_i(A)\|_2 \cdot \|\mathbf{s}_m^i\| \le C_{\rho} \left( \frac{C_2}{\rho} \right)^i,
\end{equation}
where $C_{\rho} = C_3 \tilde{C}_{\rho}$.

Under the assumption $r > C_2$, there always exists an intermediate radius $\rho$ such that $C_2 < \rho < r$. Consequently, the ratio satisfies $C_2 / \rho < 1$, which mathematically guarantees that the remainder term decays geometrically to zero:
\begin{equation*}
\lim_{i \to \infty} \left| \mathbf{u}^\top p_i(A) \mathbf{s}_m^i \right| = 0.
\end{equation*}
This rigorously completes the proof of Theorem \ref{Theorem 1}. 

\bibliographystyle{siamplain}
\bibliography{references}
\end{document}